\documentclass[12pt]{article}
\usepackage{amsfonts}
\usepackage{amsthm}
\usepackage{graphicx}
\usepackage{latexsym}
\usepackage{amssymb}
\usepackage{amsmath}
\newtheorem{definition}{Definition}[section]
\newtheorem{theorem}{Theorem}[section]
\newtheorem{lemma}{Lemma}[section]
\newtheorem{corollary}{Corollary}[section]
\newtheorem{proposition}{Proposition}[section]
\newtheorem{example}{Example}[section]
\newtheorem{remark}{Remark}[section]
\makeatletter
\begin{document}

\begin{center}
   {\bf On Predictive Density Estimation for Location Families under Integrated $L_2$ and $L_1$ Losses}
\end{center}

\begin{center}
{\sc Tatsuya Kubokawa$^{a}$, \'Eric Marchand$^{b}$, William E. Strawderman$^{c}$} \\

{\it a Department of Economics, University of Tokyo, 
7-3-1 Hongo, Bunkyo-ku, Tokyo 113-0033, JAPAN\quad 
(e-mail: tatsuya@e.u-tokyo.ac.jp) } \\
    
{\it b  Universit\'e de
    Sherbrooke, D\'epartement de math\'ematiques, Sherbrooke Qc,
    CANADA, J1K 2R1 \quad (e-mail: eric.marchand@usherbrooke.ca) } \\

{\it c  Rutgers University, Department of Statistics and Biostatistics, 501 Hill
Center, Busch Campus, Piscataway, N.J., USA, 08855 \quad (e-mail:
straw@stat.rutgers.edu) }
\end{center}
\vspace*{0.2cm}
\begin{center}
{\sc Summary} \\
\end{center}
\small
Our investigation concerns the estimation of predictive densities and a study of efficiency  as measured by the frequentist risk of such predictive densities with integrated $L_2$ and $L_1$ losses.  Our findings relate to a $p-$variate spherically symmetric observable $X \sim p_X(\|x-\mu\|^2)$ and the objective of estimating the density of $Y \sim q_Y(\|y-\mu\|^2)$ based on $X$.  For $L_2$ loss, we describe Bayes estimation, minimum risk equivariant estimation (MRE), and  minimax estimation.  We focus on the risk performance of the benchmark minimum risk equivariant estimator, plug-in estimators, and plug-in type estimators with expanded scale.
For the multivariate normal case, we make use of a duality result with a point estimation problem bringing into play reflected normal loss.   In three of more dimensions (i.e., $p \geq 3$), we show that the MRE estimator is inadmissible under $L_2$ loss and provide dominating estimators.  This brings into play Stein-type results for estimating a multivariate normal mean with a loss which is a concave and increasing function of $\|\hat{\mu}-\mu\|^2$.  We also study the phenomenon of improvement on the plug-in density estimator of the form $q_Y(\|y-aX\|^2)\,,
0<a \leq 1\,,$ by a subclass of scale expansions $\frac{1}{c^p} \, q_Y(\|(y -aX)/c \|^2)$
with $c>1$, showing in some cases, inevitably for large enough $p$, that all choices $c>1$ are dominating estimators.
Extensions are obtained for scale mixture of normals including a general inadmissibility result of the MRE estimator for $p \geq 3$.   Finally, we describe and expand on analogous plug-in dominance results for spherically symmetric distributions with 
$p \geq 4$ under $L_1$ loss.   

\vspace*{0.5cm}
\small
\noindent  {\it AMS 2010 subject classifications:}   62C20, 62C86, 62F10, 62F15, 62F30

\noindent {\it Keywords and phrases}:  Concave loss, Convolutions, Dominance, Frequentist risk, Inadmissibility, L1 loss, L2 loss, Minimax, Minimum risk equivariant, Multivariate normal, Predictive density, Restricted parameter space, Scale mixture of normals, Stein estimation.

\normalsize
\section{Introduction}

\subsection{The model and problem}

\noindent Consider independently distributed
\begin{equation}
\label{model}
X|\mu \sim p(x-\mu), Y|\mu \sim q(y-\mu)\,; x,y,\mu \in \mathbb{R}^p;
\end{equation}
where $p$ and $q$ are known, not necessarily equal, and $\mu$ is unknown.  
\noindent  For predictive analysis purposes, researchers
are interested in specifying a predictive density $\hat{q}(y|x)$ as an estimate of the density
$q(y-\mu)$.  In turn, such a density may play a surrogate role for generating either future or missing values of $Y$.   \\
\noindent Our interest and motivation here lies in assessing the efficiency of such predictive densities with integrated $L_2$
and $L_1$ losses and corresponding frequentist risk, where
\begin{eqnarray}
\nonumber
L_{\alpha}(\mu,\hat{q}) = \int_{\mathbb{R}^p} |q(y-\mu) - \hat{q}(y)|^{\,\alpha} \, dy\,, \\
R_{\alpha}(\mu,\hat{q}) = \int_{\mathbb{R}^p} \{ \int_{\mathbb{R}^p} |q(y-\mu) - \hat{q}(y)|^{\,\alpha} dy \}\, p(x-\mu) \, dx\,,
\end{eqnarray}
for $\alpha=1,2$.   The set-up in (\ref{model}) includes the normal model with
\begin{equation}
\label{normalmodel}
X|\mu \sim N_p(\mu, \sigma^2_X I_p), Y|\mu \sim N_p(\mu, \sigma^2_Y I_p)\,,
\end{equation}
scale mixtures of normal distributions (Definition \ref{definition1}),  
and more generally spherically symmetric distributions with
\begin{equation}
\label{ss}
X|\mu \sim p_X(\|x-\mu\|^2)\,,\,Y|\mu \sim q_Y(\|y-\mu\|^2)\,,
\end{equation}
to which the developments of this paper will relate.

\subsection{Motivation and overview of findings}

\noindent Our research work is motivated by the need for understanding structural elements of this problem, and, as expanded upon below, our findings focus mainly on: {\bf (A)} the benchmark minimum risk equivariant (MRE) estimator $\hat{q}_{mre}$, {\bf (B)} the performance of plug-in estimators $q(y-\hat{\mu}(X))$ of $q(y-\mu)$, $y \in \mathbb{R}^p$, where $\hat{\mu}(X)$ is an estimator of $\mu$, and {\bf (C)} improvements on plug-in estimators obtained by expanding the scale (or variance). 
 
\begin{enumerate}
\item[ {\bf (A)}]  In our problem, the MRE estimator for $L_2$ loss is obtained by the generalized Bayes estimator of the density $q(y-\mu)$ with respect to the flat prior 
$\pi(\mu)=1$ on $\mathbb{R}^p$.  Furthermore, it is minimax and thus represents an important benchmark and an attractive choice as an estimator.   These features also hold for Kullback-Leibler (KL) loss (e.g., Liang and Barron, 2004; Kubokawa et al., 2013), defined as $L_{KL}(\mu, \hat{q})=\int_{\mathbb{R}^p}  q(y-\mu) \,\log(\frac{q(y-\mu)}{ \hat{q}(y)}) \,dy$.  Although $\hat{q}_{mre}$ can possess other interesting features, such as being an admissible estimator for normal models (\ref{normalmodel}) under KL loss and with $p=1,2$ (Brown, George, Xu, 2008), Komaki (2001) established the inadmissibility of $\hat{q}_{mre}$ for such a normal model, KL loss, for $p \geq 3$, as well as provide dominating estimators.  With a striking parallel between this result and Stein's inadmissibility of the sample mean as a point estimator of the mean $\mu$ of a $N_p(\mu, \sigma_X^2 I_p)$ population under squared error loss, further analytical relationships between these predictive density estimation and point estimation problems were obtained by 
George, Liang and Xu (2006), Brown, George and Xu (2008), and Fourdrinier et al. (2010) among others.   \\

\noindent  In Section \ref{normalsection}, we arrive at a replication of Komaki's inadmissibility of 
$\hat{q}_{mre}$ result under $L_2$ loss for the normal model (\ref{normalmodel}) with $p \geq 3$, and provide dominating estimators.  We further extend the result to scale mixtures of normals in Section \ref{mresmnsection}.  These results are achieved by first establishing key relationships between our predictive density estimation problem and a  problem of estimating $\mu$ based on $X \sim p(x-\mu)$ under a loss of the type $f(\|\hat{\mu}-\mu\|^2)$ where $f$ (which depends on $q$) is shown to be increasing and 
concave.  Then, we capitalize on known results and/or familiar techniques (e.g., Brandwein and Strawderman, 1991, 1981; Brandwein, Ralescu, Strawderman, 1993) for obtaining dominating point estimators of the usual procedure $X$, which thus lead to dominating predictive density estimators of $\hat{q}_{mre}$, and the latter's inadmissibility.  The dual loss functions that intervene are of independent interest on their own and our findings also represent contributions from the point estimation perspective.  Namely, the dual loss for the normal model turns out to be reflected normal loss (Spiring, 1993).  Finally, we also provide, for $L_2$ loss, various properties and examples relative to $\hat{q}_{mre}$.  These properties also apply for KL loss since the MRE estimators coincide and a key representation of $\hat{q}_{mre}$ involves a convolution of $p$ and $q$ in (\ref{model}) (Proposition \ref{mreprop}, Example \ref{bayes} (a)).

\item[ {\bf (B)}]  Plug-in estimators are ubiquitous in statistical theory and practice.
For the univariate normal model (\ref{normalmodel}) and KL loss, Aitchison (1975) showed  that the flat prior Bayes procedure (which is $\hat{q}_{mre}$) is a $N(x, \sigma^2_X + \sigma^2_Y)$ density, and furthermore showed that it dominates the plug-in $N(x, \sigma^2_Y)$ density.  Lawless and Fredette (2005) present an instructive approach using a pivotal quantity to obtain KL improvements on plug-in estimators as well. Fourdrinier et al. (2010) elaborate on plug-in estimators $q(y-\hat{\mu}(X))$ for normal models and KL loss.  Their inadmissibility may be directly attributable, in some cases, to the inadmissibility of $\hat{\mu}(X)$ in estimating $\mu$ under (a dual) squared error loss
(also see part {\bf C} below for another factor explaining their inefficiency under KL loss).    \\

\noindent  In Section \ref{L1section}, we study the performance of plug-in estimators for $L_1$ loss and spherically symmetric models.  By making use of a key identity (Lemma \ref{pluginL1}), we obtain a dual point estimation loss.  We pursue this with inferences concerning the plug-in
density $q(t-X)\,, \, t \in \mathbb{R}^p$, establish its inadmissibility quite generally for $p \geq 4$, and obtain dominating plug-in estimators of the form $q(t-\hat{\mu}(X))$.  This is achieved in a similar manner as in Section \ref{mresmnsection} and as described in {\bf (A)} by using the fact that the dual loss is an increasing and concave function of squared error $\|\hat{\mu}-\mu \|^2$.  A larger class of dominating estimators is obtained in Section \ref{L1smnsection} for scale mixtures of normals.  
For $L_2$ loss, we do not deal as explicitly with plug-in estimators of the form $q(y-X)$ since these are invariant and are thus dominated by the MRE estimator $\hat{q}_{mre}$ under $L_2$ loss.  This explains our focus in {\bf (A)} on rather providing improvements of 
$\hat{q}_{mre}$ (if possible) under $L_2$ loss.

\item[ {\bf (C)}]  Fourdrinier et al. (2010) show, for normal model plug-in estimators $q(y-\hat{\mu}(X))$ and KL loss, that a range of scale expansions (or variance expansions) always lead to improvements of the form $\hat{q}_c(y;X)=\frac{1}{c^p} q(\frac{y-\hat{\mu}(X)}{c})$ with $c>1$ regardless of the plugged-in estimator $\hat{\mu}(X)$ and the dimension $p$ as long as it is not degenerate.  This may appear paradoxical since the variance associated with the plug-in density $q(y-\hat{\mu}(X))$ matches the variance of the true density $q(y-\mu)$, but it is always best to ignore this true variance and to opt for an estimator $\hat{q}_c$ whose associated variance overestimates the true variance.   From the loss function perspective, this is also somewhat paradoxical in that the estimate $\hat{\mu}(x)$ approaches $\mu$, the loss associated with the plug-in $\hat{q}_1$ approaches $0$, while the losses associated with other $\hat{q}_c$'s do not approach $0$.

\noindent We obtain various findings extending this phenomenon to $L_2$ loss with $\hat{\mu}(X)=aX$: for normal models and $a=1$ (Section \ref{duality}), 
normal models and $0<a <1$ (Section \ref{aXsection}), scale mixtures of normal distributions and $a=1$ (Section \ref{scaleexpansionsection}).  In Section \ref{duality},  the unbiased predictive density estimator, which is of the form $\hat{q}_c(y;X)$ is also improved on.  A surprise arises : in some cases, typically when the dimension $p$ is large enough, all expansions $\hat{q}_c(y;X)$ with $c>1$ improve on the plug-in estimator $\hat{q}_1(y;X)$ !    As an example, for normal cases with equal variances, this unusual situation occurs for all $p \geq 4$.  Taking $c$ to be infinitely large is of course silly as it becomes equivalent to using a flat density estimate converging to $0$, but the $L_2$ penalty is bounded in the normal case (and in some generality), however silly your estimate, and the result brings home another point of view on the inefficiency of the plug-in estimator.
  
\end{enumerate}

\noindent   Other findings (Corollaries \ref{ms2005}, \ref{Hartigan}, \ref{ms2005L1}) in this paper relate directly to restricted parameter space settings, where $\mu$ belongs to some known subset of $\mathbb{R}^p$, and 
are derived by exploiting dual relationships between predictive density and point estimation problems as well as restricted parameter space findings (e.g., Marchand and Strawderman, 2004).  The paper is organized as follows.   Section \ref{defsection} contains definitions and properties relative to convolutions, scale mixtures of normals, and the $L_2$ distance between two multivariate normal densities.  The latter technical result is extended to a spherically symmetric setting (Section \ref{extension}) and to $L_1$ loss (Section \ref{L1section}).  Sections \ref{bayes+mresection} and \ref{minimaxsection} focus on Bayes, best equivariant, and minimax estimation, with properties and accompanying examples.
The developments of Sections \ref{duality} and \ref{aXsection} relate to themes {\bf (A), (B), (C)} described above and to the multivariate normal model (\ref{normalmodel}) and $L_2$ loss.  Section \ref{extension} extends several results of Section \ref{normalsection} from multivariate normal to scale mixtures of multivariate normal models, including a $p \geq 3$ inadmissibility result for $\hat{q}_{mre}$ (Section \ref{mresmnsection}) and improvements by expansion of scale (Section \ref{scaleexpansionsection}).   Section \ref{L1section} deals with $L_1$ loss and is highlighted by a $p \geq 4$ inadmissibility result for the plug-in density $\hat{q}(y-X)$ and the specifications of dominating predictive density estimators given for general spherically symmetric cases (Section \ref{L1identitysection}), as well as for scale mixtures of normals (including normal) (Section \ref{L1smnsection}).  Several of the approaches taken in Sections \ref{normalsection}, \ref{extension} and \ref{L1section}  are analogous and involve dual multivariate location vector point estimation problems under various loss functions $l(\|\hat{\mu}-\mu\|^2)$ with $l$ generated from $q$, and with $l$ increasing and concave.  
Although our primary applications relate to the predictive density estimation problem, several of our findings represent point estimation findings, complement existing results, are of interest on their own.

\section{Definitions, preliminary results.  Bayes, best equivariant and minimax estimation under $L_2$ loss}
\label{prelsection}
\subsection{Some definitions and preliminary results}
\label{defsection}
\noindent We collect here various definitions and properties which will be useful throughout the paper.
Distributions in (\ref{model}) include the subclass of scale mixture of normals, with examples given by the multivariate Cauchy, Student, Logistic, Laplace, Generalized Hyperbolic, Exponential Power distributions, among others (e.g., Andrews and Mallows, 1974).  

\begin{definition}
\label{definition1}
Model (\ref{model}) is referred to as a scale mixture of normals model whenever
\begin{equation}
\label{sn}  p(x) = \int_{\mathbb{R}_+} \phi(\frac{x}{v^{1/2}})\, v^{-p/2}\, dG(v)\,,\, 
q(y) = \int_{\mathbb{R}_+} \phi(\frac{y}{w^{1/2}})\, w^{-p/2}\, dH(w)\,,
\end{equation}
for $x,y \in \mathbb{R}^p$,  where $\phi$ is (hereafter) taken to be the normal $N_p(0,I_p)$ density, and 
$W \sim G$, $V \sim H$ are independently distributed mixing random variables on 
$\mathbb{R_+}$, for which we further assume that $E(V^{-p/2})$ and $E(W^{-p/2})$ are finite.  We will denote such models or densities as $p \sim SN_p(G)$ and $q \sim SN_p(H)$.  
\end{definition}
\noindent  Convolutions $p*q$ will be omnipresent in this paper (e.g., Lemma \ref{repr})
and are given by $p*q(t) \,=\, \int_{\mathbb{R}^p} \, p(t-u) \, q(u) \, du $, $t \in\mathbb{R}^p$, for densities $p$ and $q$.   
Just as it is the case for the subclass of normal distributions, the above subclass of scale mixture of normals is closed with respect to convolutions.  

\begin{lemma}
\label{convolutionscalemixture}
For $p \sim SN_p(G)$ and $q \sim SN_p(H)$, we have  $p*q \sim SN_p(F)$ where $F$ is the cumulative distribution function (cdf) of $T=^d V+W$.
\end{lemma}
\noindent
{\bf Proof.}  Since, conditionally on $(V,W)$, $X$ and $Y$ are independently distributed 
as $N_p(0,V I_p)$  $N_p(0,W I_p)$ respectively, it follows that $X+Y|V,W \sim N_p(0,(V+W) I_p)$ whence the result.  \qed

\noindent The following result, of which the latter part gives the $L_2$ distance between multivariate normal densities, will be used several times.  A generalization  is given below in Lemma \ref{L2distance}.
\begin{lemma}  
We have for all $\mu_1, \mu_2 \in \mathbb{R}^p$ and $\sigma_1, \sigma_2 \in \mathbb{R}_+$:
\begin{equation}
\label{identity}
\int_{\mathbb{R}^p} \phi(\frac{y-\mu_1}{\sigma_1}) \, \phi(\frac{y-\mu_2}{\sigma_2}) \, dy \, =\, (\frac{\sigma_1^2 \, \sigma_2^2}{\sigma_1^2+\sigma_2^2})^{p/2} \, \phi (\frac{\mu_1-\mu_2}{\sqrt{\sigma_1^2+\sigma_2^2}} ) \,, 
\end{equation}
\begin{equation}
\label{normalL2distance}
 \, \int_{\mathbb{R}^p}   \, (\frac{1}{\sigma_1^p} \, \phi(\frac{y-\mu_1}{\sigma_1}) \,- \frac{1}{\sigma_2^p} \, \phi(\frac{y-\mu_2}{\sigma_2})\,
)^2  \, dy  = \frac{1}{(4\pi \sigma_1^2)^{p/2}} +\frac{1}{(4\pi \sigma_2^2)^{p/2}} \,-\, \frac{2}{(\sigma_1^2 + \sigma_2^2)^{p/2}}\, \phi (\frac{\mu_1-\mu_2}{\sqrt{\sigma_1^2+\sigma_2^2}} )\,.
\end{equation}
\end{lemma}
\noindent {\bf Proof.}  Identity (\ref{identity}) is readily verifed.  For (\ref{normalL2distance}), expand the square on the left-hand side to obtain 
\begin{equation*}
\, \frac{1}{ (\sigma_1^2)^{p}} \, \int_{\mathbb{R}^p} \phi^2(\frac{y-\mu_1}{\sigma_1}) \, dy \, + \, \frac{1}{(\sigma_2^2)^{p}} \int_{\mathbb{R}^p} \phi^2(\frac{y -\mu_2}{\sigma_2})    \, dy
- \frac{2}{(\sigma_1 \sigma_2)^{p}} \int_{\mathbb{R}^p} \, \phi(\frac{y-\mu_1}{\sigma_1}) \, \phi(\frac{y -\mu_2}{\sigma_2}) \, dy\,.
\end{equation*}  Applying identity (\ref{identity}) to these three terms leads to  (\ref{normalL2distance}).  \qed 

\subsection{Bayes and minimum risk equivariant estimators}
\label{bayes+mresection}
\noindent As in the case of Kullback-Leibler loss, Bayes estimators under $L_2$ loss are simply given by the predictive density
$q(y|x)$.  

\begin{lemma}
\label{bayes}
For model (\ref{model}), integrated $L_2$ loss, a prior density $\pi$ for $\mu$, and a posterior density $\pi(\mu|x)$ with respect to measure $\nu$, the Bayes predictive density estimator of $g(y-\mu)$, $y \in \mathbb{R}^p$, is given by
\begin{equation} 
\label{bayesformula}
\hat{q}_{\pi}(y;x)= \int_{\mathbb{R}^p} q(y-\mu) \, \pi(\mu|x) \, d\nu(\mu)\,.
\end{equation}
\end{lemma}
\noindent {\bf Proof.}  The posterior loss for estimator $\hat{q}(\cdot)$ is given by
\begin{equation}
\nonumber
\int_{\mathbb{R}^p} \{ \int_{\mathbb{R}^p} (q(y-\mu) - \hat{q}(y))^{\,2} dy \}\, \pi(\mu|x) \, d\nu(\mu)\,.
\end{equation}
\noindent Interchanging the order of integration, we see that for each $y$ the minimizing $\hat{q}(y)$ is the posterior expectation
$E^{\mu|x} (q(y-\mu))$ which, being a density as a function of $y$, yields the result.  \qed 

\noindent  For location models as in (\ref{model}) with spherically symmetric $q$, we obtain an interesting representation when the posterior density is location invariant.
\begin{lemma}
\label{repr}
In Lemma \ref{bayes}, whenever the posterior density is location invariant of the form $\pi(\mu|x)=g(\mu-\hat{\mu}(x))$, the Bayes predictive density estimator of $q(y-\mu)$, $y \in \mathbb{R}^p$, is equal to $q*g(y-\hat{\mu}(x))$, where $q*g$ is the convolution of $q$ and $g$.  
\end{lemma}
\noindent {\bf Proof.}
From (\ref{bayesformula}) and with $\pi(\mu|x)=g(\mu-\hat{\mu}(x))$, the Bayes predictive density of $q(y-\mu)$ is equal to
$\int_{\mathbb{R}^p} q(y-\mu) \,g(\mu-\hat{\mu}(x)) \, d \mu\, = \int_{\mathbb{R}^p} q(y-\hat{\mu}(x)- \mu') \, g(\mu') \, d\mu' = (q*g)(y-\hat{\mu}(x)).$  \qed

\begin{remark}
Since the Bayes predictive density estimators coincide for Kullback-Leibler and $L_2$ losses, the above lemma and the examples presented at the end of this subsection apply as well to Kullback-Leibler loss.
\end{remark}
\noindent The minimum risk equivariant estimator  can be derived as the Bayes rule with respect to the Haar invariant prior $\pi(\mu)=1$ for $\mu$, and is minimax.  This follows as the problem is invariant under the group of location changes (including the choice of loss), and from a general representation for the minimum risk equivariant estimator as the Bayes estimator associated with the corresponding Haar measure (e.g., Eaton, 1989), and with the minimaxity following from Kiefer (1959).  The following Proposition summarizes the above and provides a direct and  instructive approach in deriving the minimum risk equivariant (mre) estimator under $L_2$ loss, which is analogous to results obtained by Murray (1977) or Kubokawa et al. (2013) for Kullback-Leibler loss.

\begin{proposition} 
\label{mreprop}
The minimum risk equivariant estimator of $q(y-\mu)$, $y \in \mathbb{R}^p$, for model (\ref{model}) and $L_2$ loss is given by
\begin{equation}
\label{mreformula}
\hat{q}_{\hbox{mre}}(y;x) \,=\, q*\bar{p} (y-x)\,,
\end{equation}
with $\bar{p}(t)=p(-t)$ for all $t$, and matches the Bayes predictive density with respect to the uniform prior on $\mathbb{R}^p$ given in Example \ref{bayesexamples} {\bf (a)}.
Furthermore, $\hat{q}_{\hbox{mre}}(\cdot;X)$ is a minimax predictive density estimator.
\end{proposition}
\noindent {\bf Proof.}  While Kiefer's result, mentioned above, gives minimaxity quite 
generally for the MRE estimator, minimaxity is established directly in Section \ref{Appendix} for the
general location case via the argument of Girschick and Savage (1951), 
using a (least favourable) sequence of Uniform priors
on the product sets  $\{\mu: |\mu_i| < k/2\,, i=1, \ldots, p\}$, $k=1,2, \ldots $ . We give a similar direct, but simpler argument specifically  for the normal case in Section 2.3.

For the minimum risk equivariance property, we only need to establish (\ref{mreformula}).  First, equivariant estimators under the additive group of transformation $(x,y) \to (x+a,y+a)$; $a \in \mathbb{R}^p$; satisfy the identity
\begin{equation}
\label{form}
\hat{q}(y;0) = \hat{q}(y-x;0) \hbox{ for all } x,y \in \mathbb{R}^p\,,
\end{equation}
 as seen by setting $a=-x$.  The risk of such estimators is constant in $\mu \in \mathbb{R}^p$ and given by  
\begin{eqnarray*}
R(\mu,\hat{q}) &=& \int_{\mathbb{R}^p} \{\int_{\mathbb{R}^p} (q(y-\mu) - \hat{q}(y-x;0))^2 dy\} \; p(x-\mu) dx \\
\, &=& \int_{\mathbb{R}^p} \{\int_{\mathbb{R}^p} (q(y) - \hat{q}(y-x;0))^2 dy\} \; p(x) dx\, \\
\, &=& \int_{\mathbb{R}^p} \{\int_{\mathbb{R}^p} (q(u+v) - \hat{q}(v;0))^2 \,p(u) \,du \} \; dv\,,
\end{eqnarray*}
with the last equality obtained with transformation $(x,y) \to (u=x,v=y-x)$. 
Now, for all $v \in \mathbb{R}^p$, the inner integral above is minimized by choosing $\hat{q}(v;0)$ to be the expected value
of $q(v+U)$ with $U \sim p$, i.e. $\hat{q}_{\hbox{opt}}(v)=\int_{\mathbb{R}^p} \, q(u+v) \, p(u) \, du\,.$  Finally, this along with (\ref{form})
tell us that
$$ \hat{q}_{\hbox{mre}}(y|\,x) = \hat{q}_{\hbox{opt}}(y-x) =   \int_{\mathbb{R}^p} \, q(u+y-x) \, p(u) \, du\, = \int_{\mathbb{R}^p} \, q(y-x-u) \, p(-u) \, du\,=\,
q*\bar{p}(y-x)\,. \qed$$
\begin{example}
\label{bayesexamples}
\begin{enumerate}
\item[ {\bf (a)}] Consider model (\ref{model}) with the uniform prior $\pi(\mu)=1$ on $\mathbb{R}^p$ and with the corresponding Bayes predictive density estimator coinciding with the MRE estimator (see Proposition \ref{mreprop}).  This gives us: $x-\mu|x \sim p$ and Lemma \ref{repr} applies with $g(y)=\bar{p}(y)=p(-y)$ and the representation
$\hat{q}_{\hbox{mre}}(y;x)= (q*\bar{p})(y-x)$.  Moreover, if $p$ is spherically symmetric as in (\ref{ss}), we have $\hat{q}_{\hbox{mre}}(y;x)= (q*p)(y-x)$.  

\item[ {\bf (b)}]
For the normal case (\ref{normalmodel}) with $q \sim N_p(0,\sigma_Y^2 I_p)$, $p \sim N_p(0,\sigma_X^2 I_p)$, we obtain immediately that $q*p \sim N_p(0, (\sigma_X^2 +\sigma_Y^2) I_p)$ and that $\hat{q}_{\hbox{mre}}(y;x)= N_p(x, (
\sigma_X^2 +\sigma_Y^2) I_p)$. 

\item[ {\bf (c)}]  Now, consider the normal case (\ref{normalmodel}) with $\pi(\mu) \sim N_p(\theta,\tau^2 I_p)$.  
Since $\mu|x \sim N_p(\hat{\mu}(x), (\tau')^2 I_p)$ with $\hat{\mu}(x) =
\frac{\tau^2 x}{\sigma_X^2 +\tau^2} + \frac{\sigma_X^2 \theta}{\sigma_X^2 +\tau^2}$ and $(\tau')^2 = \frac{\sigma_X^2 \tau^2}{\sigma_X^2 +\tau^2}$, this corresponds in the notation of Lemma \ref{repr} to $q \sim N_p(0,\sigma_Y^2 I_p)$, $g \sim N_p(0, (\tau')^2 I_p)$.  The same Lemma \ref{repr} thus implies that $\hat{q}_{\pi}(y;x) \sim N_p(\hat{\mu}(x), (\sigma_Y^2 + (\tau')^2) I_p)$.  As seen from above, we further point out that all predictive densities $N_p(aX+b, (\sigma_Y^2 + a \sigma_X^2) I_p)$ with $0 \leq a <1$, $b \in \mathbb{R}$ are unique Bayes estimators with finite Bayes risks and hence admissible.  We will show in Section \ref{normalsection} that the MRE estimator (i.e., $a=1, b=0$) is inadmissible for $p \geq 3$. 

\item[ {\bf (d)}]  As a further illustration of {\bf (a)} and extension of {\bf (b)}, consider the uniform prior $\pi(\mu)=1$ and scale mixtures of normals densities $p \sim SN_p(G)$, $q \sim SN_p(H)$ as in (\ref{sn}).  It thus follows from part {\bf (a)} of this Example and Lemma \ref{convolutionscalemixture} that 
\begin{equation}
\hat{q}_{\hbox{mre}}(y;x) = \int_{\mathbb{R}_+} \phi(\frac{y-x}{t^{1/2}})\, t^{-p/2}\, dF(t)\,,
\end{equation}
with $F$ the cdf of $T=^d W+V$.

\item[ {\bf (e)}]  (Multivariate Student and Cauchy models)  A prominent scale mixture of normals example is the multivariate Student $T(\nu,\sigma)$ with degrees of freedom $\nu>0$ and scale parameter $\sigma>0$.  In (\ref{model}), this corresponds to 
$X-\mu \sim p \sim T(\nu_1,\sigma_1)$,  $Y-\mu \sim q \sim T(\nu_2,\sigma_2)$ where the density of  $T(\nu,\sigma)$ is given by
\begin{equation*}
 \frac{\Gamma(\frac{1}{2}(p+\nu))}{(\pi \nu \sigma^2)^{p/2}\,\,\Gamma(\frac{\nu}{2})} \; (1 + \frac{\|t\|^2}{\nu \, \sigma^2})^{-\frac{1}{2}(p+\nu))}\,.
\end{equation*}
Part {\bf (a)} tells us that $\hat{q}_{\hbox{mre}}(y;x)= (q*p)(y-x)$.  Such a convolution density, including cases where one of the densities is that of a normal distribution,  has arisen in other settings and been analyzed by others (e.g., Nason 2006; Berg and Vignat, 2010).   
The particular case of a multivariate Cauchy ($\nu_1=\nu_2=1$) gives rise, simply, to
\begin{equation*}
\hat{q}_{\hbox{mre}}(y;x) \,=\, \frac{\Gamma(\frac{1}{2}(p+1))}{\pi \,\,(\sigma_1 + \sigma_2)^p} \; (1 + \frac{\|y-x\|^2}{(\sigma_1 + \sigma_2)^2})^{-\frac{1}{2}(p+1)}\,
\end{equation*}
since $T(1,\sigma_1) * T(1, \sigma_2) = T(1, \sigma_1+\sigma_2)$. 

\item[ {\bf (f)}]  (Exponential location model-  Kubokawa et al., 2013.) As a further illustration of Lemma \ref{repr}, or of part {\bf (a)} of this sequence of examples, consider $X=(X_1, \ldots, X_n)'$, where the $X_i$'s and $Y$ are independently distributed as $\hbox{Exp}(\mu, \beta_1)$, $Y \sim \hbox{Exp}(\mu, \beta_2)$, with densities $\beta_i^{-1} e^{-\frac{(t-\mu)}{\beta_i}} \, \mathbb{I}_{(\mu,\infty)}(t)$ and known $\beta_1, \beta_2$.  For the uniform prior $\pi(\mu)=1$, we obtain as the posterior density $\pi(\mu|x) = \frac{1}{n\beta_1} e^{\frac{\mu -\hat{\mu}(x)}{n\beta_1}}\, \mathbb{I}_{(-\infty, 0)}(\mu -\hat{\mu}(x))$, with $\hat{\mu}(x)=\min_i x_i$.
Lemma \ref{repr} thus applies with $g(u)= \frac{1}{n\beta_1} e^{\frac{u}{n\beta_1}}\, \, \mathbb{I}_{(-\infty, 0)}(u )$ and $q(u) = \frac{1}{n\beta_2} e^{-\frac{u}{n\beta_2}}\, \, \mathbb{I}_{(0, \infty)}(u )$.   The convolution of $g$ and $q$ yields 
$$(q*g)(u) \, = \,  \frac{1}{n\beta_1+\beta_2} \; \{e^{-|u|/n\beta_1} \, \mathbb{I}_{(-\infty, 0)}(u ) \, + \, e^{-|u|/\beta_2}    \, \mathbb{I}_{(0, \infty)}(u ) \}\,,$$
and $$\hat{q}_{\hbox{mre}}(y;x)= (q*g)(y-\hat{\mu}(x))\,.$$
 
\item[ {\bf (g)}] (Uniform model)   Consider $X_1, \ldots, X_n, Y$ independently and uniformly distributed on the interval $(\mu, \mu+1)$.  
A uniform prior $\pi(\mu)=1$ leads to the posterior $\mu|x \sim \hbox{Uniform}(x_{(n)}-1, x_{(1)})$, with $x_{(1)}=\min_i x_i$ and 
$x_{(n)}=\max_i x_i$.  Applying (\ref{bayesformula}) directly yields
\begin{eqnarray*}
\hat{q}_{\hbox{mre}}(y;x) &=&  \int_{\mathbb{R}} \mathbb{I}_{(0,1)} (y-\mu) \, \mathbb{I}_{(x_{n}-1, x_{(1)})}(\mu) \, d\mu \\
\, &=&  (y+1-x_{(n)}) \mathbb{I}_{(x_{n}-1, x_{1}]}(y) + (x_{(1)}-x_{(n)}+1)  \mathbb{I}_{(x_{1}, x_{n}]}(y) \\
&+&  (x_{(1)} +1 -y) \mathbb{I}_{(x_{}, x_{1}+1]}(y)\,,
\end{eqnarray*}
which is a trapezoidal shaped predictive density for $Y$.   
\end{enumerate}
\end{example}

\subsection{Minimax estimator and least favourable sequence in the normal case}
\label{minimaxsection}
We provide here for normal case (\ref{normalmodel}) a direct approach to obtain a least favorable sequence of priors and show that the best equivariant estimator
$\hat{q}_{\hbox{mre}}(\cdot|X) \sim N_p(X, (\sigma_X^2 + \sigma_Y^2) \, I_p)$ is minimax under $L_2$ loss.  We proceed in a familiar way showing that 
$\hat{q}_{\hbox{mre}}(\cdot|X)$ is extended Bayes with constant risk.  We make use of the risk computation below in (\ref{riskmre}), which shows that the constant risk of 
$\hat{q}_{\hbox{mre}}(\cdot|X)$ is given by $R_0= (4\pi \sigma_Y^2)^{-p/2} \, - \, (4\pi (\sigma_X^2 + \sigma_Y^2))^{-p/2} \, .$
Consider the sequence of priors $\pi_m \sim N_p(0,m I_p)$; $m=1,2, \ldots$; as in part {\bf (c)} of Example \ref{bayesexamples} with $\theta=0, \tau^2=m$, and 
with corresponding Bayes estimators $\hat{q}_{\pi_m}(\cdot|X) \sim N_p(\hat{\mu}_{\pi_m}(X), \, \sigma_{1,m}^2\, I_p)$, $\hat{\mu}_{\pi_m}(x) = \frac{mx}{m+\sigma_X^2}$, and $\sigma_{1,m}^2 = \frac{m \sigma_X^2}{m+\sigma_X^2} + \sigma_Y^2$.  The posterior loss  $\int_{\mathbb{R}^p} (\hat{q}_{\pi_m}(y|x) - q(y-\mu )\,)^2 \,dy$ is obtained from  (\ref{normalL2distance}) with $\sigma_1^2= \sigma_{1,m}^2$, $\sigma_2^2=\sigma_Y^2$, $\mu_1=\hat{\mu}_{\pi_m}(x)$, $\mu_2=\mu$, and given by
\begin{equation*}
k_m(\mu,x) =  \frac{1}{(4\pi \sigma_{1,m}^2)^{p/2}} + \frac{1}{(4\pi \sigma_{Y}^2)^{p/2}} - \frac{2}{(\sigma_{1,m}^2 +\sigma_Y^2)^{p/2}} \, \,\phi\left(\frac{\hat{\mu}_{\pi_m}(x)-\mu}
{\sqrt{\sigma_{1,m}^2 +\sigma_Y^2}}\right) \,.
\end{equation*}
From this and by making use of (\ref{identity}), the expected posterior loss is evaluated as 
\begin{eqnarray}
\nonumber  \int_{\mathbb{R}^p} k_m(\mu)\, \pi(\mu|x)  \,d\mu &=&
\frac{1}{(4 \pi \sigma_{1,m}^2)^{p/2}} + \frac{1}{(4\pi \sigma_{Y}^2)^{p/2}} \\
\nonumber &-& \frac{2}{(\sigma_{1,m}^4 -\sigma_Y^4)^{p/2}} \, \int_{\mathbb{R}^p} \phi\left(\frac{\hat{\mu}_{\pi_m}(x)-\mu}
{\sqrt{\sigma_{1,m}^2 +\sigma_Y^2}}\right) \, \phi\left(\frac{\hat{\mu}_{\pi_m}(x)-\mu}
{\sqrt{\sigma_{1,m}^2 -\sigma_Y^2}}\right) \, dy \\
\label{bayesrisk}&=& \frac{1}{(4\pi \sigma_{1,m}^2)^{p/2}} + \frac{1}{(4\pi \sigma_{Y}^2)^{p/2}} - \frac{2}{(4\pi \sigma_{1,m}^2 )^{p/2}}\,.
\end{eqnarray}
Observe that the expected posterior loss, given in (\ref{bayesrisk}), is independent of $x$ and thus matches the Bayes risk $r_{\pi_m}$.  Finally, since 
$\sigma_{1,m}^2 \to \sigma_X^2 + \sigma_Y^2$ when $m \to \infty$, we have $\lim_{m \to \infty} r_{\pi_m} = R_0 = R(\mu, \hat{q}_{\hbox{mre}})$, which implies that
the estimator $\hat{q}_{\hbox{mre}}$ is indeed minimax (and that the sequence $\pi_m$ is least favorable).

\section{Plug-in type estimators and $L_2$ loss: the normal case}
\label{normalsection}
\subsection{Duality and the efficiency of density estimators $N_p(\hat{\mu}(X), c^2 \sigma_Y^2\,I_p)$}
\label{duality}
We consider here normal model (\ref{normalmodel}) and the performance of density estimators $\hat{q}_{c^2,\hat{\mu}} \sim N_p(\hat{\mu}(X), c^2 \sigma_Y^2\,I_p)$, which combine both a plug-in 
component with $\hat{\mu}(X)$ being an estimate of $\mu$, and a modification of variance component for $c^2 \neq 1$.
As for Kullback-Leibler loss (Fourdrinier et al. 2011), we demonstrate that the efficiency of such estimators relates to: {\bf (i)} the efficiency of the point estimator
$\hat{\mu}(X)$ in estimating $\mu$, as well as {\bf (ii)} the degree of variance expansion governed by the choice of $c^2>1$.  With respect to {\bf (i)} and the duality with the point estimation problem, it is a reflected normal loss that arises, which we denote and define as 
\begin{equation}
\label{rnl}
L_{\gamma}(\mu, \hat{\mu}) \, =\, 1 - e^{-\frac{\|\hat{\mu}(x)- \mu)\|^2}{2 \gamma}}\;, \hbox{ with } \gamma >0\,,
\end{equation}
in contrast to squared-error loss which intervenes in duality for Kullback-Leibler loss.  

\noindent  
\begin{lemma}
\label{cplugin}
For estimating a multivariate normal density of $Y \sim N_p(\mu, \sigma_Y^2 I_p)$,  the integrated $L_2$ loss
of the density estimate $\hat{q}_{c,\hat{\mu}} \sim N_p(\hat{\mu}(x), c^2 \sigma_Y^2 I_p)$ is given by
\begin{equation}
\label{cpluginformula}
\frac{1}{(\sigma_Y^2)^{p/2}} \left( \frac{1}{(4\pi)^{p/2} } + \frac{1}{(4\pi c^2)^{p/2}} -  \frac{2}{(2\pi (c^2+ 1))^{p/2}} \,  e^{-\frac{(\|\hat{\mu}(x)- \mu\|)^2}
{2 \sigma_Y^2(c^2+1)}} \, \right) ,
\end{equation}
\end{lemma}
\noindent {\bf Proof.}   This is a direct application of (\ref{normalL2distance}) with $\mu_1=\mu, \mu_2=\hat{\mu}(x)$, $\sigma_1^2=\sigma_Y^2$, and 
$\sigma_2^2=c^2 \sigma_Y^2$.  \qed  

\begin{corollary}
\label{cor}
\begin{enumerate}
\item[ {\bf (a)}]  For fixed $c^2$, the frequentist risk of the density estimator $\hat{q}_{c^2,\hat{\mu}} \sim N_p(\hat{\mu}(x), c^2 \sigma_Y^2I_p)$ under integrated $L_2$ loss
is equal to the frequentist risk of the point estimator  $\hat{\mu}(X)$ of $\mu$ under loss  $a+ b L_{\gamma_0}(\mu,\hat{\mu})$, with $a=\frac{1+ c^{p}}{(4\pi c^2 \sigma_Y^2)^{p/2}}, b= \frac{2}{(2\pi (c^2+1) \sigma_Y^2)^{p/2}}$, and $\gamma_0=(c^2+1)\,\sigma_Y^2$.  Namely, $\hat{q}_{c^2,\hat{\mu}_1} \sim N_p(\hat{\mu}_1(X), c^2 \, \sigma_Y^2I_p)$ improves on $\hat{q}_{c^2,\hat{\mu}_2} \sim N_p(\hat{\mu}_2(X), c^2 \, \sigma_Y^2I_p)$ iff $\hat{\mu}_1(X)$  improves on  $\hat{\mu}_2(X)$ under reflected normal loss $L_{\gamma_0}(\mu,\hat{\mu})$.

\item[ {\bf (b)}]  For $\hat{\mu}(X)=X$, the risk $R(\mu,\hat{q}_{c^2,\hat{\mu}})$ is constant as a function of $\mu$, and given by
\begin{equation}
\label{riskofqhatc}
\frac{1}{(2\pi \sigma_Y^2)^{p/2}} \left( (\frac{1}{2})^{p/2} \,+ \, \frac{1}{(2c^2)^{p/2}} - \frac{2}{(r+c^2+1)^{p/2}} \right)\,,
\end{equation}
with $r=\frac{\sigma_X^2}{\sigma_Y^2}$.  For all $p$, the constant (and minimax) risk of $\hat{q}_{\hbox{mre}}$, corresponding to the optimal choice $c^2=1+r$ is equal to 
\begin{equation}
\label{riskmre} R(\mu, \hat{q}_{\hbox{mre}}) \, =\, \frac{1}{(4\pi \sigma_Y^2)^{p/2}} \, - \, \frac{1}{(4\pi (\sigma_X^2 + \sigma_Y^2))^{p/2}} \, .  
\end{equation} 
\item[ {\bf (c)}]
Furthermore, all estimators $\hat{q}_{c^2,\hat{\mu}}$ with $c^2>1$ dominate the plug-in or mle estimator $\hat{q}_{1,\hat{\mu}}$ whenever $p \geq p_0=\frac{\log 4}{\log (1+ r/2)}$, and, otherwise for $p < p_0$, the estimator $\hat{q}_{c,\hat{\mu}}$  dominates $\hat{q}_{1,\hat{\mu}}$ iff $1 < c^2 < k(p,r)$, where $k(p,r)$ is the unique solution on $(1,\infty)$ in $c^2$ of the equation
\begin{equation}
\label{root}
(\frac{1}{2})^{p/2} + 2 (\frac{1}{r+c^2+1})^{p/2} - 2 (\frac{1}{r+2})^{p/2} - (\frac{1}{2c^2})^{p/2}\,=\,0\,.
\end{equation}
\end{enumerate}
\end{corollary}
\noindent {\bf Proof.}   Part {\bf (a)} follows directly from Lemma \ref{cplugin}.  For part {\bf (b)}, use (\ref{cpluginformula}) with
$\hat{\mu}(x)=x$ and $\hat{q}_{c^2} \sim N_p(x,c^2\,I_p)$ to obtain
\begin{equation*}
R(\mu, \hat{q}_{c^2}) =  \frac{1}{(\sigma_y^2)^{p/2}} \left(\frac{1}{(4\pi)^{p/2}} + \frac{1}{(4\pi c^2)^{p/2}} - 
\frac{2}{(c^2+1)^{p/2}} \, E_{\mu}\left(\phi(\frac{X-\mu}{\sqrt{\sigma_Y^2(c^2+1)}})  \right) \right)\,.
\end{equation*}
The use of identity (\ref{identity}) leads to (\ref{riskofqhatc}).  Now, set 
$\psi(c^2) = \frac{1}{(2c^2)^{p/2}} - \frac{2}{(r+c^2+1)^{p/2}}$
so that $(2\pi \sigma_y^2)^{p/2} \, R(\mu,\hat{q}_{c^2}) = (\frac{1}{2})^{p/2} + \psi(c^2) $, and observe that $\hbox{sgn} (\psi'(c^2)) = \hbox{sgn} \{(2c^2)^{1+p/2} - (r+c^2+1)^{1+p/2}\}$ for $c>0$.  From this, we infer that $\psi'(c^2)$ changes signs once, on $(0,\infty)$, from $-$ to $+$ at $c^2=1+r$, which along with the evaluation of
(\ref{riskofqhatc}) for $c^2=1+r$ establishes part {\bf (b)}.  For part {\bf (c)}, the comparison of $\hat{q}_{c^2}$ with $\hat{q}_1$ for $c^2 > 1$ hedges on the sign of $\psi(c^2)-\psi(1)$ for $c > 1$.  From above, we know that $\psi(c^2) -\psi(1)$ is either negative for all $c^2 >1$, or negative iff $1< c^2 < k(p,r)$.  Finally, we have $\lim_{c^2 \to \infty} \{\psi(c^2)-\psi(1)\} = \frac{2}{(2+r)^{p/2}} - \frac{1}{2^{p/2}} \leq 0$ if and only if $p \geq p_0$, concluding the proof.  \qed
 
\begin{example}
\label{r=1}
For equal variances (i.e., $r=1$) $\sigma_X^2$ and $\sigma_Y^2$,  we obtain $p_0 \approx 3.419$ so that universal dominance for all choices $\hat{q}_{c^2,\hat{\mu}}$ with $c^2>1$ over $\hat{q}_{1,\hat{\mu}}$ arises for $p \geq 4$.  
And for $p \geq 3$, the cut-off points $k(1,p)$ are given by
$k(1,2)=6$ (exact), $k(1,1) \approx 4.65$, and $k(1,3) \approx 11.47$. We remark upon the fact that $p_0$ decreases as the ratio $r=\frac{\sigma_X^2}{\sigma_Y^2}$ increases so that the above universal dominance occurs also for (at least) all $p \geq 4$ whenever $\sigma_X^2 > \sigma_Y^2$.  For instance if $\sigma_X^2=2\sigma_Y^2$ (i.e., $r=2$), we obtain $p_0=2$.  

\end{example}

\noindent The following is a consequence of part {\bf (b)} of Corollary \ref{cor}.
\begin{remark}
\label{mre}
For integrated $L_2$ loss and model (\ref{normalmodel}), the ratio of risks between the minimum risk equivariant estimator $\hat{q}_{1+r,X} \sim N_p(X, (\sigma_X^2 + \sigma_Y^2) \, I_p)$ and the plug-in estimator $\hat{q}_{1,X} \sim N_p(X, \sigma^2_Y I_p)$ is given by $\frac{R(\mu,\hat{q}_{1,X})}{R(\mu,\hat{q}_{1+r,X})}\,=\,2 \left(\frac{1 - (\frac{1}{1+r/2})^{p/2} }{1  -  (\frac{1}{1+r})^{p/2} }  \right)$.  It is easy to verify that this ratio increases in both $r$ and $p$, and approaches $2$ when either $r$ or $p$ increase to $\infty$.   The monotonicity in the ratio of variances $r=\frac{\sigma_X^2}{\sigma_Y^2}$ translates, understandably, to worsening performance of $\hat{q}_1$ as the relative variability of the observable $X$ increases.  
\end{remark} 

\begin{remark} (On the unbiased predictive density estimator)
\label{umvue}
An unbiased predictive density estimator of the density of $Y|\mu$, (i.e., of $\frac{1}{(\sigma_Y^2)^{p/2}} \, \phi(\frac{y-\mu}{\sigma_Y})$) exists whenever $\sigma_X^2 < \sigma_Y^2$ (for the univariate case, see for instance Lehmann and Casella, 1998; or Shao, 1999).  Indeed, considering density estimates 
$q_{c^2,X} \sim N_p(X,  c^2 \sigma_Y^2 I_p)$, we have from (\ref{identity}) 
$$  \int_{\mathbb{R}^p}  \frac{1}{(c^2 \sigma_Y^2)^{p/2}} \, \phi(\frac{y -x}{c \, \sigma_Y}) \, \frac{1}{(\sigma_X^2)^{p/2}} \, \phi(\frac{x -\mu}{ \sigma_X})\,dx=\,
\frac{1}{(\sigma_X^2 + c^2 \sigma_Y^2)^{p/2}} \, \phi(\frac{y -\mu}{\sqrt{(\sigma_X^2 + c^2 \sigma_Y^2)}})\,,$$
so that a $N_p(X,  c^2 \sigma_Y^2 I_p)$ density is an unbiased estimator of a $N_p(\mu, \sigma_X^2 +  c^2 \sigma_Y^2 I_p)$ density (pointwise and globally), and the choice $ c^2=1-\frac{\sigma_X^2}{\sigma_Y^2} (>0)$ yields an unbiased estimator of the density of $Y|\mu$.  Since $X$ is a complete sufficient statistic, it follows that this estimator is the sole unbiased estimator.
\footnote{The more standard setup, perhaps, has $\sigma_X^2=\frac{\sigma_Y^2}{n}$, where $n$ is the size of a sample drawn from $X$.}
Here, the unbiased predictive density estimator shrinks the variance, instead of expanding it.  It will thus, with its risk given by (\ref{riskofqhatc}) and as already analysed as a function of $ c^2$, perform even worse than the plug-in $\hat{q}_{1,X}$.  In fact, it is dominated by the plug-in, the best equivariant estimator, a range of choices $\hat{q}_{c,X}$, $1-r< c^2 < k_0(p,r)$, and with $k_0(p,r)=+\infty$ as soon as $p \geq - \frac{\log 4}{\log (1-r)}$. 
\end{remark}

\subsection{Improvements over the minimum risk equivariant estimator}
\label{mreimprovementsnormalsection}
\noindent As presented in part {\bf (a)} of Corollary \ref{cor}, there is a dual point estimation loss to predictive density estimation for plug-in estimators in the sense that the frequentist risk of a plug-in density estimator is given by the frequentist point estimation risk of the same plug-in estimator under reflected normal loss $a+bL_{\gamma_0}$.  
Reflected normal loss was introduced by Spiring (1993), namely as an option for a bounded loss.  It is also not convex in $\|d-\mu\|$, but strictly bowled shaped in $\|d-\mu\|$. 
We can thus borrow results applicable to such loss functions. For instance, results from Marchand and Strawderman (2005), or again Kubokawa and Saleh (1994), show that for $p=1$ the Bayes estimator $\hat{\mu}_{\pi_U}(X)$ with respect to the uniform prior either on a compact interval $(a,b)$ or left-bounded interval $(a,\infty)$ dominates the MRE estimator $X$ under strictly bowled shaped loss and hence reflected normal loss.  Here is a formulation of such an inference.

\begin{corollary}
\label{ms2005}
 For estimating a univariate normal density $Y \sim N(\mu, \sigma_Y^2)$ based on $X \sim N(\mu, \sigma_X^2)$ under integrated $L_2$ loss, and with the restriction $\mu \in [a,b]$ ($\mu \in [a,\infty)$), the estimator  $\frac{1}{c \, \sigma_Y } \, \phi(\frac{y -\hat{\mu}_{\pi_U}(X)}{c \, \sigma_Y })$ dominates $\frac{1}{c \, \sigma_Y } \, \phi(\frac{y -X}{c \sigma_Y })$, where $\hat{\mu}_{\pi_U}(X)$ is the Bayes point estimator of $\mu$ associated with a uniform prior on $[a,b]$ (on $[a,\infty$)) under reflected normal loss $L_{\gamma_0}(\mu,\hat{\mu})$ with $\gamma_0= (c^2+1) \sigma_Y^2$.
\end{corollary}
\noindent {\bf Proof.}  Since $\hat{\mu}_{\pi_U}(X)$ dominates the MRE estimator $X$ as shown by Marchand and Strawderman (2005), the result is a consequence of 
part {\bf (a)} of Corollary \ref{cor}. \qed \\

\noindent
Another class of applications of part {\bf (a)} of Corollary \ref{cor} are generated by
estimators $\hat{u}(X)$ that dominate $X$, for $X \sim N_p(\mu, \sigma_X^2 I_p)$ and for 
loss functions $f(\|d-\mu\|^2)$ with $f$ increasing and concave.  Such findings were given 
by Brandwein and Strawderman (1991, 1981), as well as Brandwein, Ralescu, and Strawderman (1993), and apply for the above reflected normal loss.  It is interesting that such an example arises naturally here in view of the fact that Brandwein and Strawderman's (1981) examples (but not theory) concern $L_p$ loss with $0<p<2$.  The developments that follow make use of similar techniques but exploit the specific nature of the loss function to obtain a wider class of dominating estimators for $p \geq 3$ of $\hat{q}_{\hbox{mre}}$.

\begin{lemma}
\label{stein}
Let $X \sim N_p(\mu, \sigma_X^2I_p)$ with known $\sigma_X^2$, and consider estimating $\mu$ under reflected normal loss loss $L_{\gamma}(\mu,\hat{\mu})=1- e^{-\frac{\|\hat{\mu}- \mu\|^2}{2 \gamma}}\,. $
Then $\hat{\mu}(X)$ dominates $X$ under $L_{\gamma}$ whenever $\hat{\mu}(Z)$ dominates $Z$ for the model $Z \sim N_p(\mu, \frac{\gamma \sigma_X^2}{\gamma+\sigma_X^2}\,\,I_p)$ under loss
$\|\hat{\mu}- \mu\|^2$.
\end{lemma}
\noindent {\bf Proof.}  Since 
$- e^{-\frac{\|\hat{\mu}- \mu\|^2}{2\gamma}} = - e^{-\frac{\|x- \mu\|^2}{2\gamma}}\, \times \, 
e^{-\frac{1}{2\gamma}  (\|\hat{\mu}- \mu\|^2 - \|x- \mu\|^2)}$, it is seen that
$$L_{\gamma}(\mu,\hat{\mu}(x)) = 1 - e^{-\frac{\|x- \mu\|^2}{2\gamma}} + e^{-\frac{\|x- \mu\|^2}{2\gamma}} 
\left(1- e^{-\frac{1}{2\gamma}  (\|\hat{\mu}(x)- \mu\|^2 - \|x- \mu\|^2)}  \right)\,,$$
for all $x \in \mathbb{R}^p$.  In terms of the risk $R_{\gamma}(\mu,\hat{\mu})=E_{\mu}(L_{\gamma}(\mu,\hat{\mu}(X))$, we thus have
$$R_{\gamma}(\mu,\hat{\mu}) = R_{\gamma}(\mu,X) + E^X_{\mu} \left[e^{-\frac{\|X- \mu\|^2}{2\gamma}} 
\left(1- e^{-\frac{1}{2\gamma}  (\|\hat{\mu}(X)- \mu\|^2 - \|X- \mu\|^2)}  \right)\right]\,,$$
and
$$\Delta_{\gamma}(\mu) = R_{\gamma}(\mu,\hat{\mu}) - R_{\gamma}(\mu,X) = (\frac{\gamma}{\gamma+\sigma_X^2})^{p/2} \, 
E^Z_{\mu} \left(1- e^{-\frac{1}{2\gamma}  (\|\hat{\mu}(Z)- \mu\|^2 - \|Z- \mu\|^2)}  \right) \,,$$
where $Z \sim N_p(\mu, \frac{\gamma\sigma_X^2}{\gamma+\sigma_X^2}\,I_p)$.  Now, note that $1-e^{-d} < d$ for any $d  \neq 0$, so that
\begin{equation}
\label{boundDelta}
\Delta_{\gamma}(\mu) < (\frac{\gamma}{\gamma+\sigma_X^2})^{p/2} \, E^Z_{\mu}\left(\|\hat{\mu}(Z)- \mu\|^2 - \|Z- \mu\|^2)  \right)\,,
\end{equation}
which yields the result.  \qed  \\

\noindent Corollary \ref{cor}'s duality between the performance {\bf (I)} of plug-in density estimators under $L_2$ loss and
the performance {\bf (II)} of the corresponding point estimator under reflected normal loss, coupled with the previous lemma which links the 
latter's point estimation performance {\bf (II)} with the one under squared error loss {\bf (III)},  lead to the following inadmissibility result and comparisons
for our predictive density estimation problem {\bf (I)}.

\begin{corollary} 
\label{mreimprovements}
For estimating a multivariate normal density $Y \sim N_p(\mu, \sigma_Y^2I_p)$ based on $X \sim N_p(\mu, \sigma_X^2I_p)$ under integrated $L_2$ loss, the estimator $\hat{q}_{\hbox{mre}}(\cdot;X) \sim N_p(X, (\sigma_X^2 + \sigma_Y^2) I_p)$ is inadmissible for $p \geq 3$, and dominated 
by any $\hat{q}(\cdot;X) \sim N_p(\hat{\mu}(X),(\sigma_X^2 + \sigma_Y^2 ) I_p)$, where $\hat{\mu}(Z)$ dominates $Z$ for $Z \sim N_p(\mu, \sigma_Z^2 \, I_p)$ under loss $\|\hat{\mu}- \mu\|^2$ and with $\sigma_Z^2 = \frac{(2\sigma_Y^2 + \sigma_X^2) \sigma_X^2}{2(\sigma_Y^2+\sigma_X^2)}$.
\end{corollary}
\noindent {\bf Proof.}  This is a direct consequence of part {\bf (a)} Corollary \ref{cor} and Lemma \ref{stein},
 applied for $c^2=1+\frac{\sigma_X^2}{\sigma_Y^2}$ and $\gamma=(c^2+1) \sigma_Y^2 \,=\, 2 \sigma_Y^2 + \sigma_X^2$.  \qed

\noindent   The above results establishes the inadmissibility of $\hat{q}_{\hbox{mre}}$ for $p \geq 3$.  Along with Stein estimation findings under squared error loss, we can generate explicit dominating plug-in type densities $\hat{p} \sim N_p(\hat{\mu}, (\sigma_X^2 + \sigma_Y^2)I_p)$.  Here are some examples.

\begin{example}
If $\pi$ is  superharmonic prior, the Bayes estimator $\hat{\mu}_{\pi}(Z)$ dominates $Z$ for $Z \sim N_p(\mu, \sigma_Z^2 I_p)$ for $p \geq 3$ as an estimator of $\mu$ under squared error loss (Stein 1981), and the corresponding plug-in density $N_p(\hat{\mu}_{\pi}(X), (\sigma_X^2 + \sigma_Y^2) I_p)$ dominates $\hat{q}_{\hbox{mre}}$ under $L_2$ loss.  In terms of the secondary problem of estimating $\mu$ under reflected normal loss $L_{\gamma}$, Lemma \ref{stein} implies the dominance of $\hat{\mu}_{\pi}(X)$ over $X$.  Wider classes of dominating estimators arise, for instance, by requiring that $\sqrt{m_{\pi}(z)}$ be superharmonic where $m_{\pi}(z)$ is the marginal density of $Z$ under $\pi$ (Fourdrinier, Strawderman, Wells, 1998).  
\end{example}

\begin{example}  
\label{baranchikexample}
Another class of dominating estimators of $Z \sim N_p(\mu, \sigma_Z^2 I_p)$ for $p \geq 3$, which will arise Section \ref{extension} as well, are given by 
Baranchik type estimators (Baranchik, 1971) $\hat{\mu}_{a,r(\cdot)}(Z) = (1 - a \frac{r(Z'Z)}{Z'Z}) \,Z$, such that $r(\cdot)$ is an increasing function, $0 \leq r(\cdot) \leq 1$, $r(\cdot) \neq 0$, and $0 < a \leq 2 (p-2) \, \sigma_Z^2\,$.  In view of Corollary \ref{mreimprovements}, the corresponding density  estimators  $N_p(\hat{\mu}_{a,r(\cdot)}(X), (\sigma_X^2 + \sigma_Y^2) I_p)$ dominate 
$\hat{q}_{\hbox{mre}}$ for $p \geq 3$ under $L_2$ loss.
\end{example}

\noindent Further applications of Corollary \ref{mreimprovements} include the following  Hartigan type result for cases where $\mu$ is restricted to $C$, $C$ being a strict subset of $\mathbb{R}^p$ which is convex with a non-empty interior.  Such cases include restrictions to balls and to cones such as order constraints $\mu_1 \leq \mu_2 \leq \ldots \leq \mu_p$, tree order constraints $\mu_i \geq \mu_1$ for all $i$, etc.

\begin{corollary}
\label{Hartigan}
Let $Z \sim N_p(\mu, \sigma_Z^2\,I_p)$ with $\sigma_Z^2= \frac{(2\sigma_Y^2 + \sigma_X^2) \sigma_X^2}{2(\sigma_Y^2+\sigma_X^2)}$, and let $\hat{\mu}_{\pi}(Z)$ be the Bayes estimator of $\mu$ associated with prior density $\pi$ and loss
$\|\hat{\mu}- \mu\|^2$.  For estimating the density of $Y \sim N_p(\mu, \sigma_Y^2 I_p)$ based on $X \sim N_p(\mu, \sigma_X^2 I_p)$ under integrated $L_2$ loss, and for the restriction $\mu \in C$ with $C$ a convex subset of $\mathbb{R}^p$ with non-empty interior:

\begin{enumerate}
\item[ {\bf (a)}]  the estimator $\hat{q}(\cdot;X) \sim N_p(\hat{\mu}_{\pi_U}(X), (\sigma_X^2 + \sigma_Y^2)\,I_p)$ dominates $\hat{q}_{\hbox{mre}}(\cdot;X) 
\sim N_p(X, (\sigma_X^2 + \sigma_Y^2)\,I_p)$ with $\pi_U$ being the uniform prior on $C$;
\item[ {\bf (b)}]  for the univariate case with $C=[a,b]$, dominance of $\hat{q}_{\hbox{mre}}(\cdot;X)$ is achieved by any $\hat{q}(\cdot;X)
\sim N_p(\hat{\mu}_{\pi}(X), (\sigma_X^2 + \sigma_Y^2)\,I_p)$
as long as the prior density $\pi$ is absolutely continuous, symmetric about $(a+b)/2$, and increasing and logconcave on $[\frac{a+b}{2},b]$.
\end{enumerate}
\end{corollary}  
\noindent {\bf Proof.}  The results are a consequence of part {\bf (a)} of Corollary \ref{cor} and Lemma \ref{stein} (with $c^2=1+r$) paired with point estimation results of
Hartigan (2004) for {\bf (a)}, and Kubokawa (2005) or Marchand and Payandeh Najafabadi (2011) for {\bf (b)}.  \qed

\begin{remark}
In the context of Corollary \ref{mreimprovements}, for cases where $\|\mu\| \leq m$, as well as cases where $\mu \in C$ with $C$ a convex cone, 
$\|\hat{\mu}_{\hbox{mle}}(x) -\mu\|^2$ is stochastically smaller than $\|x-\mu\|^2$ for all $\mu$ so that  $\hat{\mu}_{\hbox{mle}}(X)$ dominates $X$ as an estimator of
$\mu$ under loss $L_{\gamma}$, $\gamma >0$.  Therefore $\hat{q}_{c^2,\hat{\mu}_{\hbox{mle}}} \sim N_p(\hat{\mu}_{\hbox{mle}}(X), c^2 \,\sigma_Y^2 \,I_p)$
dominates $\hat{q}_{c^2,\hat{\mu}_{0}} \sim N_p(X, c^2\,\sigma_Y^2\,I_p)$ for such restricted parameter spaces, $c^2 >0$.
\end{remark}

\subsection{Plug-in estimators with $\hat{\mu}(X)=aX$:  improvements by expanding the scale}
\label{aXsection}
\noindent  We revisit here the normal case (\ref{normalmodel}) and analyze the performance of the estimators $\hat{q}_{c^2,\hat{\mu}} \sim N_p(\hat{\mu}(X),c^2 \sigma_Y^2 I_p)$ with more development for the
affine linear case $\hat{\mu}(x)=ax$, $0 < a \leq 1$.  As seen in Section \ref{duality}, there exists for $a=1$ an optimal choice (i.e., $c^2=1+r$ with $r=\frac{\sigma_X^2}{\sigma_Y^2}$) of the expansion factor $c^2$ and, for $p \geq p_0=\frac{\log 4}{\log (1+r/2)}$, we can expand the variance as much as desired and still dominate the plug-in 
$N_p(x,\sigma_Y^2 I_p)$.  The objective here is to assess whether such results hold for other choices of $\hat{\mu}(X)$ and more specifically: {\bf (i)} to determine a range of variance expansions or values $c^2$ that lead to improvement, and {\bf (ii)} to
determine whether there exists a universal dominance result for sufficiently large $p$, where universal means for all $c^2 >1$.  Explicit findings with respect to {\bf (i)} and Kullback-Leibler loss were obtained by Fourdrinier et al. (2011) with the maximum amount of allowable expansion to retain improvement for all $\mu$ an increasing function of the infimum squared error risk.    \\
 
\noindent We start off with the risk expression (see Lemma \ref{cplugin})

\begin{equation}
\label{risk}
R(\mu, \hat{q}_{c^2,\hat{\mu}}) = \frac{1}{(4\pi \sigma_Y^2)^{p/2}} + \frac{1}{(2\pi \sigma_Y^2)^{p/2}} 
\, \left(\frac{1}{(2\,c^2)^{p/2}} \, -\, \frac{2}{(c^2+1)^{p/2}} \, E_{\mu} (e^{-\frac{\|\hat{\mu}(X)-\mu  \|^2}{2 \sigma_Y^2 (c^2+1)}}) \right)\,, 
\end{equation} 
and the derivative 
\begin{equation}
\label{drisk}
\frac{\partial}{\partial c^2} \, R(\mu, \hat{q}_{c^2,\hat{\mu}}) \, = \, \frac{1}{(2\pi \sigma_Y^2)^{p/2}} \left(\frac{1}{(c^2+1)^{p/2+1}} 
\, E_{\mu}\left( e^{-\frac{Z}{2(c^2+1)}} \, (p - \frac{Z}{c^2+1})\right)  -\, \frac{p}{(2c^2)^{p/2+1}} \,,  \right) 
\end{equation} 
with $Z=\frac{\| \hat{\mu}(X)- \mu  \|^2}{\sigma_Y^2}$.

\begin{remark}
\label{alwayslocally}
Plug-in estimators $\hat{q}_{1,\hat{\mu}} \sim N_p(\hat{\mu}(X), \sigma_Y^2 I_p)$ with non-degenerate $\hat{\mu}(X)$ can always be improved locally at $\mu$ by an 
expansion $\hat{q}_{c^2,\hat{\mu}} \sim N_p(\hat{\mu}(X),c^2 \sigma_Y^2 I_p)$  with $c^2 \in (1, c_0^2(\mu))$ for some $c_0^2(\mu)>1$ (dependent on $\hat{\mu}(X)$).  If possible, global dominance
$\mu \in \Theta$ is thus achieved by selecting $c^2 \in (1, \inf_{\mu \in \Theta} c_0^2(\mu)]$.   This can be seen by a continuity argument and (\ref{drisk}), since  $$\frac{\partial}{\partial c^2} \, R(\mu, \hat{q}_{c^2,\hat{\mu}})|_{c^2=1} \, = \, \frac{1}{2} \, \frac{1}{(4\pi \sigma_Y^2)^{p/2}} 
\left( E_{\mu} (e^{-\frac{Z}{4}} (p- \frac{Z}{2}) \, -\,  p \right) < 0\,,$$
as $e^{-\frac{y}{4}} (p- \frac{y}{2}) \, -\,  p \leq 0$ for all $y \geq 0$, with equality iff $y=0$. 
\end{remark}

\noindent For the particular case $\hat{\mu}(X)=aX, 0<a<1$, we arrive at more explicit expressions for the risk and its derivative in (\ref{risk}) and (\ref{drisk}) by using the 
exact distributional result $Z=\frac{\|aX-\mu\|^2}{\sigma_Y^2} \sim a^2 \,r\, \chi^2_p (\frac{(a-1)^2  \|\mu\|^2 }{a^2 \sigma_X^2} )$, and the mixture representation :
$$ Z|L \sim Gamma (\frac{p}{2} + L, 2 a^2 \, r)\,,\, L \sim \hbox{Poisson}(\,\frac{\delta}{2}\,), \hbox{ with } \delta=\frac{(a-1)^2  \|\mu\|^2 }{a^2 \sigma_X^2}, \, r=
\frac{\sigma_X^2}{\sigma_Y^2}\,.$$

\begin{lemma}
\label{comp}
For $Z = \frac{\|aX-\mu  \|^2}{\sigma_Y^2}$, we have
\begin{equation}
E_{\mu} \left( e^{-\frac{Z}{2(c^2+1)}} \, (p - \frac{Z}{c^2+1})  \right) \,= \, (p-h) \, \theta^{\frac{p}{2}+1} \, e^{-\frac{h}{2}}\,,
\end{equation}
with $h=\frac{(a-1)^2 \, \|\mu\|^2}{a^2 \sigma_X^2+(c^2+1) \sigma_Y^2}$ and $\theta=\frac{c^2+1}{a^2 r +c^2+1}$.
\end{lemma}
\noindent {\bf Proof.}  Straightforward calculations yield  
$$ E_{\mu} \left( e^{-\frac{Z}{2(c^2+1)}} \, (p - \frac{Z}{c^2+1}) | L=l \right) \, =\, p \theta^{\frac{p}{2}+l+1} \, +\, 2l \,\theta^{\frac{p}{2}+l+1}\,
(1-\frac{1}{\theta})\,, $$ with $L \sim \hbox{Poisson}(\delta/2)$.   The result follows by the Poisson related evaluations $E( \theta^{\frac{p}{2}+L+1}) =
\theta^{\frac{p}{2} +1} \, e^{-h/2}$ and $E( L \, \theta^{\frac{p}{2}+L+1}) =
\theta^{\frac{p}{2} +1} \frac{\theta \delta}{2} \, e^{-h/2}$, and by collecting terms. \qed  \\

\noindent Setting $\psi_a(c^2)= \frac{\partial}{\partial c^2} \, R(\mu, \hat{q}_{c^2,\hat{\mu}})$ for $\hat{\mu}(X)=aX$, it thus follows from the above expression
and (\ref{drisk}) that

\begin{equation}
\label{psi}
\psi_a(c^2) \,=\, \frac{1}{(2\pi \sigma_Y^2)^{p/2}} \, \frac{1}{(a^2 r +c^2+1)^{p/2+1}} \, \left((p-h) \, e^{-h/2} \, - \, p \, (\frac{a^2 r+c^2+1}{2c^2})^{p/2+1}  \right)\,,
\end{equation}
with $h$ as in Lemma \ref{comp}.

\begin{lemma}
\label{riskvaries}
For all $a \in (0,1)$, $\mu \in \mathbb{R}^p$,  $\psi_a(c^2)$ changes signs once from $-$ to $+$ as 
$c^2$ increases on $[1,\infty)$. 
\end{lemma}
\noindent {\bf Proof.}   We have already established (Remark \ref{alwayslocally}) that $\psi_a(1)<0$, and $\psi_a(c^2)$ is clearly negative for $h \geq p$, i.e., $1 \leq c^2 \leq k_0, $ with $k_0=\left(1 \wedge (\frac{(a-1)^2 \, \|\mu\|^2}{p\sigma_Y^2} - a^2r-1))\right)$.  Otherwise, it is easy to see that $\psi_a(c^2)$ is increasing in $c^2$ for $c^2 > k_0$.  Finally, the result follows since $\lim_{c^2 \to \infty} ((a^2+c^2+1)^{p/2+1}) (2\pi \sigma_Y^2)^{p/2} \, \psi_a(c^2) \, = \, p (1- (\frac{1}{2})^{p/2+1}) >0.  $  \qed

\begin{theorem}
\label{dominance}
For $\hat{\mu}(X)=aX$, $0<a<1$, $\hat{q}_{c^2,\hat{\mu}}$ dominates $\hat{q}_{1,\hat{\mu}}$ under integrated $L_2$ loss if and only if $1 < c^2 \leq k_a(p)$, where $k_a(p)=\infty$ whenever $p \geq p_0(a)=\frac{2 \, \log(2)}{\log(1+\frac{a^2 r }{2})}$, and otherwise when $p < p_0(a)$, $k_a(p)$  is the unique solution in $c^2 \in (1,\infty)$ of the equation    
\begin{equation}
\label{roota}
(\frac{1}{2})^{p/2} + 2 (\frac{1}{c^2+a^2 r +1})^{p/2} - 2 (\frac{1}{a^2r+2})^{p/2} - (\frac{1}{2c^2})^{p/2}\,=\,0\,.
\end{equation}

\end{theorem}

\begin{remark}
\label{gen}
For $a=1$, we recover part {\bf (c)} of Corollary \ref{cor} and, namely, the universal in $c^2$ dominance for $p \geq p_0(1)$.  Observe that the universal dominance property for all $c^2>1$ is inevitable for large enough $p$.   This is not necessarily the case for other estimators (see footnote below).
\end{remark}

\noindent {\bf Proof of Theorem \ref{dominance}.} 
We accompany the proof with several observations.  Here are the elements of the proof.

\begin{enumerate}
\item[ (A)]  Locally at $\mu$, it follows from Lemma \ref{riskvaries} and for $\hat{\mu}(X)=aX$ that
\begin{equation}
\label{local}
R(\mu, \hat{q}_{1,\hat{\mu}}) - R(\mu, \hat{q}_{c^2,\hat{\mu}}) \geq 0 \; \hbox{ if and only if } 1 < c^2 \leq c_0^2(\mu)\,,
\end{equation}
with equality only at $c^2=c_0^2(\mu)$.  
\item[ (B)]  It is possible that $c_0^2(\mu)=+\infty$, not only for $\hat{\mu}(X)=aX$ but also more generally for other $\hat{\mu}$'s.  In fact, a necessary and general condition for this type of universal dominance is as follows.  

\begin{lemma}
\label{i}
For any non-degenerate $\hat{\mu}(X)$ and $T=\frac{\|\hat{\mu}(X)-\mu  \|^2}{\sigma_Y^2}$, we have 
\begin{equation}
\label{infinity}
R(\mu, \hat{q}_{1,\hat{\mu}}) - \lim_{c \to \infty} R(\mu, \hat{q}_{c^2,\hat{\mu}}) \geq 0 \hbox{ for all } \mu \; \hbox{ if and only if } \sup_{\mu} E_{\mu} (e^{-T/4}) \leq \frac{1}{2}\,.  \footnote{With Jensen's inequality and focussing at $\mu=0$, a necessary condition for this is $E_{0}(T) \geq 4 \log_e(2)$.  For $p \geq 3$ and the James-Stein estimator $\hat{\mu}_{JS}(X)= (1 -\frac{(p-2)\sigma_X^2}{X'X}) X$, we have $E_{0}(T)=\frac{1}{\sigma_Y^2}\;E_0(\|\hat{\mu}_{JS}(X)-0\|^2=2r $, so that dominance at $c\to \infty$ is not possible for $\hat{\mu}$ taken to be the James-Stein estimator, or any other estimator dominating $\hat{\mu}_{JS}(X)$ such as its positive part, whenever $r < 4 \log_e 2.$}
\end{equation}
\end{lemma}

\noindent {\bf Proof.}  From (\ref{risk}), we obtain $(4\pi \sigma_Y^2)^{p/2}\, \lim_{c^2 \to \infty} R(\mu, \hat{q}_{c^2,\hat{\mu}})=1$, and 
$(4\pi \sigma_Y^2)^{p/2}\, R(\mu, \hat{q}_{1,\hat{\mu}})=2-2\,E_{\mu} (e^{-T/4})\,, $ yielding the result.  \qed

\item[ (C)]  Applying Lemma \ref{i} for $\hat{\mu}(X)=aX$, we make use of the stochastically increasing property of the family of distributions
of $Z=\|aX-\mu\|^2 \,\sim a^2 r \chi^2_p (\frac{(a-1)^2 \|\mu\|^2}{a^2 \sigma_X^2})$, with $\|\mu\|^2$ viewed as the parameter, to infer that
$$ \sup_{\mu \in \mathbb{R}^p} E_{\mu}(e^{-Z/4})=E_0(e^{-Z/4})\,=\, (1+\frac{a^2r}{2})^{-p/2}\,.$$  Therefore, condition (\ref{infinity}) becomes equivalent to
$(1+\frac{a^2r}{2})^{-p/2} \leq \frac{1}{2}$ which is $p \geq p_0(a)$.  

\item[ (D)]  We set $\Delta_{c^2}(\mu)= R(\mu, \hat{q}_{1,\hat{\mu}}) - R(\mu, \hat{q}_{c^2,\hat{\mu}})$.  Using Karlin's variation diminishing properties (e.g., Brown, Johnstone and MacGibbon, 1981) and a monotone likelihood property of the family of non-central Chi-square distributions, we obtain the following global comparison.

\begin{lemma}
\label{global}  For $\hat{\mu}(X)=aX$, $a \in (0,1)$, we have for all fixed $c^2>1$:  $\Delta_{c^2}(\mu) \geq 0$ for all $\mu \in \mathbb{R}^p$ if and only if
$\Delta_{c^2}(0) \geq 0$.  Furthermore, $\Delta_{c^2}(0) \geq 0$ if and only if $1 <c^2 \leq k_a(p)$ (as stated in Theorem \ref{dominance}).
\end{lemma} 
\noindent {\bf Proof.}  First, the condition $\Delta_{c^2}(0) \geq 0$ is equivalent to $1 <c^2 \leq k_a(p)$ as can be seen by making use of the risk expression in (\ref{risk}) and the evaluation $E_0(e^{-Z/d})=(1+\frac{2a^2 r}{d})^{-p/2}$ (here for $\mu=0$, $Z \sim \hbox{Gamma}(p/2,2a^2 r)$) applied for $d=4$ and $d=2(c^2+1)$.  Secondly, we obtain  from (\ref{risk}) $(4\pi \sigma_Y^2)^{p/2} \, \Delta_{c^2}(\mu) \,=\, E_{\|\mu\|^2} (g(Z))$, with $$g(z)= 1 - c^{-p} + 2 \left( (\frac{2}{c^2+1})^{p/2}\, e^{-\frac{z}{2(c^2+1)}}  - e^{-\frac{z}{4}} \right) \,.$$
With (i) $g'(z) \geq 0$ if and only if $z \leq \frac{(p+2) \, \log(\frac{c^2+1}{2})}{\frac{1}{2} - \frac{1}{c^2+1}}$, (ii) $g(0)=  2 (\frac{2}{c^2+1})^{p/2}- 1 - c^{-p} <0$ for all $p,c^2 \geq 1$, and (iii)  $\lim_{z \to \infty} g(z) = 1 - c^{-p} >0$ for all $p \geq 1, c^2>1$, we infer that 
$g(z)$ changes signs exactly once from $-$ to $+$ as $z$ increases on $\mathbb{R}_+$.  Consequently, the variation diminishing properties applicable to the family of distributions of $Z$, which has an increasing monotone likelihood property in $Z$ with parameter $\|\mu\|^2$, imply that $E_{\|\mu\|^2} (g(Z))$ changes signs at most once as a function 
of $\|\mu\|^2 \in \mathbb{R}_+$ whence the result.  \qed

With the steps above, we have now established Theorem \ref{dominance}.   \qed

\begin{remark}  For $\mu=0$, the optimal value ${c^*}^2$ of $c^2$ is available from (\ref{psi}) and the equation  $\psi({c^*}^2)=0$ yielding ${c^*}^2=1+a^2r$.
We thus have $k_a(p) \geq 1+a^2r$ in view  of the previous Lemma.  
\end{remark} 
\end{enumerate}

\section{Extensions to scale mixtures of normals under integrated $L_2$ loss}
\label{extension}
The developments in this section parallel those of Section \ref{normalsection} but relate to scale mixtures of normals.   In Section \ref{L2identiysection}, which applies more generally for multivariate location families, we obtain an explicit representation for the $L_2$ distance (Lemma \ref{L2distance}) between two densities of the same multivariate location family which brings into play a convolution related to these densities. 
This leads to a duality between the $L_2$ risk of the MRE estimator, and more generally for density estimators of the form  $f(y-\hat{\mu}(X))$, $y \in \mathbb{R}^p$,  of a density $q(y-\mu)$ with a risk function (under a certain loss) of $\hat{\mu}(X)$ as a point estimator of $\mu$.
For many spherically symmetric choices of $q$, namely all scale mixtures of normals, the dual point estimation resulting loss is an non-decreasing and concave function of $\|\hat{\mu}-\mu\|^2$ and, as in Section \ref{mreimprovementsnormalsection}, we establish in Section \ref{mresmnsection}, for $p \geq 3$ and with risk finiteness conditions, the inadmissibility of $\hat{q}_{\hbox{mre}}$ as well as provide dominating estimators.   Finally, in Section \ref{scaleexpansionsection}, we assess the risk performance of scale expansion estimators of the form $\frac{1}{c^p}\, q(\frac{y-x}{c})$, $c>1$, in comparison with the plug-in estimator $q(y-x)$, and replicate some of the normal case features with improvements always to be found in this subclass.

\subsection{An identity for $L_2$ distance and general dominance results of plug-in type predictive density estimators}
\label{L2identiysection}
We begin this section with a general $L_2$ distance identity conveniently expressed in terms of convolutions.

\begin{lemma}
\label{L2distance}
Whenever finite, the $L_2$ distance  $\rho_{L_2}=\int_{\mathbb{R}^p} |\, q(t-\mu_1) - f(t-\mu_2)|^2 \,dt \,$ between densities $f(t-\mu_1)$ and $q(t-\mu_2)$, $ \mu_1, \mu_2 \in \mathbb{R}^p$, is given by $\rho_{f,q}(\mu_2-\mu_1)$ with
\begin{equation}
\label{L2edistanceformula}
\rho_{f,q}(s)=   q*\bar{q}\,(0) + f*\bar{f}\,(0)- 2 q * \bar{f}(s), \; s \in \mathbb{R}^p\,,
\end{equation}
$\bar{q}(t)=q(-t)$, and $\bar{f}(t)=f(-t)$ for all $t \in \mathbb{R}^p$. 
\end{lemma}
\noindent {\bf Proof.}  In a straightforward manner, we have 
\begin{eqnarray*}
\rho_{L_2}
\,&=& \int_{\mathbb{R}^p} q^2(t-\mu_1) \, dt + \int_{\mathbb{R}^p} f^2(t-\mu_2) \, dt - 2 \int_{\mathbb{R}^p} q(t-\mu_1) \, f(t-\mu_2) \, dt \\
\, &=& \int_{\mathbb{R}^p} q(-t)\,\bar{q}(t) +  \int_{\mathbb{R}^p} f(-t)\,\bar{f}(t) \, dt \, - \, 2 \int_{\mathbb{R}^p} q(\mu_2-\mu_1-t) \, \bar{f}(t) \, dt \\
\, &=& \,   q*\bar{q}\,(0) + f*\bar{f}\,(0)- 2 q * \bar{f}(\mu_2-\mu_1)\,\\
\, &=& \rho_{f,q}(\mu_2-\mu_1).  \qed
\end{eqnarray*}

\noindent In our predictive density estimation context, we will be seeking to estimate the density $q(y-\mu)$ under $L_2$ loss and the above provides the loss associated with the subclass of estimators of the form $f(y-\hat{\mu})$ with $f$ fixed.  Comparisons with the MRE estimator, which we carried out in the previous section for the normal case, are of particular interest.  As shown in Example 2.1, such a choice corresponds to $f \equiv q*\bar{p}$ and $\hat{\mu}(x)=x$, with $X \sim p(t-\mu)$ and $\bar{p}(t)=p(-t)$ for all $t \in \mathbb{R}^p$.   As a direct consequence of the above Lemma, we have the following Corollary which relates to the MRE estimator.

\begin{corollary}
\label{pluginL2}
For estimating the density $q(y-\mu)$, $y,\mu \in \mathbb{R}^p$, under integrated $L_2$ loss and based on $X \sim p(x-\mu)$, 
\begin{enumerate}
\item[ {\bf (a)}]  The frequentist risk of the estimator $f(y-\hat{\mu}(X))$ is equal to the frequentist risk of the point estimator $\hat{\mu}(X)$ of $\mu$ under loss 
$\rho_{f,q}(\hat{\mu}-\mu) $;
\item[ {\bf (b)}] The estimator $f(y-\hat{\mu}_1(X))$ dominates the estimator
 $f(y-\hat{\mu}_2(X))$ if and only if $\hat{\mu}_1(X)$ dominates $\hat{\mu}_2(X)$ as a point estimator of $\mu$ under loss $\rho_{f,q}(\hat{\mu}-\mu) $ or, equivalently, under loss 
\begin{equation}
\label{duallossL2}
1 - \frac{2 \, q*\bar{f}(\hat{\mu}-\mu)} {q*\bar{q}\,(0) + f*\bar{f}\,(0)}\;;
\end{equation}
\item[ {\bf (c)}] The estimator $q*\bar{p}(y-\hat{\mu}(X))$ dominates the MRE estimator  $q*\bar{p}(y-X)$ if and only if $\hat{\mu}(X)$ dominates $X$ under loss $\rho_{q*\bar{p},q}(\hat{\mu}-\mu) $ or, equivalently, under loss 
\begin{equation}
\label{duallossmreL1}
1 - \frac{2 \, q*\bar{q}*p (\hat{\mu}-\mu)}{q*\bar{q}(0) + q*\bar{p}*\bar{q}*p(0)}\,.
\end{equation}
\end{enumerate}
\end{corollary}
\noindent {\bf Proof.}  Parts {\bf (a)} and {\bf (b)} follow directly from Lemma \ref{L2distance}. 
Part {\bf (c)} follows from Proposition \ref{mreprop}'s representation of the MRE estimator and by applying part {\bf (b)} for $f=q*\bar{p}$ and $\bar{f} = \bar{q}*p$.   \qed  \\

\noindent 
\begin{remark} With the above results, Corollary \ref{cplugin} turns out to be a particular case of Lemma \ref{L2distance} by taking $q$ ($=\bar{q}$ by symmetry of $q$) $\sim 
N_p(0, \sigma_Y^2 I_p)$, $f (= \bar{¸f}) \sim N_p(0, c^2 \sigma_Y^2 I_p)$, $\mu_1=\mu$, and
$\mu_2 = \hat{\mu}$.  Indeed with these values, (\ref{L2edistanceformula}) yields expression (\ref{cpluginformula}) with the normal convolutions: $q * \bar{q} \sim N_p(0,\, 2\sigma_Y^2 I_p)$, $f * \bar{f} \sim N_p(0,\, 2 \,c^2 \,\sigma_Y^2 I_p)$, and $q * \bar{f} \sim N_p(0, \, (1+c^2) \, \sigma_Y^2 I_p)$.  Similarly, part {\bf (a)} of Corollary \ref{cor}, as well as its dual reflected normal loss  $L_{(c^2+1) \sigma_Y^2}$, follow from part {\bf (b)} of Corollary \ref{pluginL2} by setting $f \sim N_p(0,\,c^2 \sigma_Y^2I_p)$ and $q \sim N_p(0,\ \sigma_Y^2 I_p)$.   Namely, the implications of Corollary \ref{cor} relative to the MRE estimator $\hat{q}_{\hbox{mre}}(\cdot|X)$ (i.e., $c^2=1+r$, $\hat{\mu}_2(X)=X$) are given in part {\bf (c)} of Corollary \ref{pluginL2}.
\end{remark}

\subsection{Inadmissibility of the MRE density estimator for $p \geq 3$ and dominating estimators}
\label{mresmnsection}
Despite the fact that Lemma \ref{L2distance} and parts {\bf (a)} and {\bf (b)} of Corollary \ref{pluginL2} apply for general densities $q,f,p$, we will focus on new applications in part {\bf (c)} of Corollary \ref{pluginL2} for scale mixtures of normals.  As in Section \ref{mreimprovementsnormalsection}, we exploit the fact that the dual loss in (\ref{duallossL2}) is an increasing and concave function of $\|\hat{\mu}-\mu\|^2$ for scale mixtures of normals, and that dominating estimators of $X$ for such a dual loss can be derived using familiar techniques first put forth by Brandwein and Strawderman (1991, 1981), as well as Brandwein, Ralescu and Strawderman (1993). 

\noindent  The first part of the following result is an adaptation of part {\bf (c)} of Corollary \ref{pluginL2} for scale mixtures of normals and for comparing estimators with the MRE predictive density estimator $\hat{q}_{\hbox{mre}}$, the middle part establishes that point estimation dominance results with squared-error loss under an associated scale mixture of normals model generates dominating estimators of $\hat{q}_{\hbox{mre}}$, and the last part  capitalizes on an existing result (Strawderman, 1974) for scale mixtures of normals and leads to an inadmissibility result for $\hat{q}_{\hbox{mre}}$ and $p \geq 3$.

\begin{theorem}
\label{smnl2loss}
Consider estimating a scale mixture of normals density $q(y-\mu)$, $y \in \mathbb{R}^p$, with $q \sim SN_p(H)$ under $L_2$ loss and based on $X \sim p(x-\mu)$, $p \sim SN_p(G)$.
Let $W_1,W_2 \sim H$, $V_1 \sim G$ be independently distributed, and let $F$ and $J$ be the cdfs of $V_1+W_1$ and $V_1+W_1+W_2$ respectively. 
\begin{enumerate}
\item[ {\bf (a)} ]
 The estimator $q*p(y-\hat{\mu}(X))$ dominates the MRE estimator  $q*p(y-X)$, with $q*p \sim SN_p(F)$, if and only if $\hat{\mu}(X)$ dominates $X$ under loss
\begin{equation}
\label{reflectedsmnloss}
f(\|\hat{\mu} - \mu  \|^2)\,=\, K - \int_{\mathbb{R_+}} \, (2\pi \,t)^{-p/2}  
e^{- \frac{\|\hat{\mu} - \mu \|^2}{2t}}\,  \, dJ(t)\,,
\end{equation}
$K$ being a constant.
\item[ {\bf (b)} ]   The estimator $q*p(y-\hat{\mu}(X))$ dominates the MRE estimator  $q*p(y-X)$, with $q*p \sim SN_p(F)$, whenever $\hat{\mu}(X')$ dominates $X'$ under squared-error loss $\|\hat{\mu}-\mu\|^2$ and for $X'-\mu \sim SN_p(F_Z)$, $F_Z$ being the cdf of
\begin{equation}
\label{Z}
Z =^d \frac{Z_1 Z_2}{Z_1+Z_2}\,, \hbox{ with } (Z_1,Z_2) \sim d\tau(z_1,z_2) \propto \frac{1}{z_2\, (z_1+z_2)^{p/2}} \, dG(z_1) \, dJ(z_2)\,.
\end{equation}
\item[ {\bf (c)} ]  Assuming $E(Z)$ and $E(Z^{-1})$ exist,  $\hat{q}_{\hbox{mre}}$ is inadmissible for $p \geq 3$ and dominated by  $q*p(y-\hat{\mu}_{a,r(\cdot)}(X))$  with $q*p \sim SN_p(F)$, and with a Baranchik type estimator $\hat{\mu}_{a,r(\cdot)}(X) = (1 - a \frac{r(X'X)}{X'X}) \,X$ such that $r(\cdot)$ is an increasing function, $\frac{r(t)}{t}$ decreases in $t$, $0 \leq r(\cdot) \leq 1$, $r(\cdot) \neq 0$, and $0 < a \leq \frac{2 (p-2)}{ E(Z^{-1})}\,$.
\end{enumerate}
\end{theorem}
\noindent {\bf Proof.}  Part {\bf (a)}  follows from part {\bf (a)} of Corollary \ref{pluginL2} and Lemma \ref{convolutionscalemixture}'s convolution properties for scale mixtures of normals.   For part {\bf (b)}, we seek a condition that suffices for the difference in risks $\Delta(\hat{\mu},\mu) =
E_{\mu} \left[\, f (\|\hat{\mu}(X) - \mu \|^2) \, - \, f (\|X - \mu \|^2)   \,\right] $ to be less than $0$, where $f$ is given in (\ref{reflectedsmnloss}).  Since $f$ is  strictly concave, the inequality $f(s) - f(t) < f'(t) (s-t)$, for $s \neq t$ and for such $f$'s, implies for the difference in losses that
\begin{equation}
\label{differencelosses}
f(\|\hat{\mu}(x)-\mu \|^2) \, - \,  f (\|x - \mu \|^2) <
f'(\|x - \mu \|^2) \, \left(\|\hat{\mu}(x)-\mu \|^2 - \|x-\mu \|^2  \right)\,,
\end{equation}
for all $x, \mu \in \mathbb{R}^p$ such that $x \neq \hat{\mu}(x)$. 
Now, using the above, part {\bf(b)} follows since
\begin{eqnarray*}
\Delta (\hat{\mu}, \mu) &<&  \frac{1}{2} \, E_{\mu}^X  \left[ \left(\|\hat{\mu}(X)-\mu \|^2 - \|X-\mu \|^2  \right) \,  \int_0^{\infty} (2\pi t)^{-p/2} \, 
e^{-\frac{\|X -\mu \|^2}{2t}} \, \frac{dJ(t)}{t} \right] \\
\, &=&  \frac{1}{2} \, \int_{\mathbb{R}^p} \left(\|\hat{\mu}(x)-\mu \|^2 - \|x-\mu \|^2  \right)  \,\int_0^{\infty} \int_0^{\infty} 
(4\pi^2 \,st)^{-\frac{p}{2}}
 \, e^{-\frac{\|x -\mu \|^2}{\frac{2st}{s+t}}}  \, dG(s) \, \frac{dJ(t)}{t} \, dx \\
& \propto &  \, \int_{\mathbb{R}^p} \left(\|\hat{\mu}(x)-\mu \|^2 - \|x-\mu \|^2  \right)  \,\int_0^{\infty} \int_0^{\infty} 
( \,\frac{2\pi\,st}{s+t})^{-\frac{p}{2}} \,
 \, e^{-\frac{\|x -\mu \|^2}{\frac{2st}{s+t}}}  \, d\tau(s,t) \, dx  \\
& \propto &  \, \int_{\mathbb{R}^p} \left(\|\hat{\mu}(x)-\mu \|^2 - \|x-\mu \|^2  \right)  \,\int_0^{\infty}
(2\pi z)^{-\frac{p}{2}}
 \, e^{-\frac{\|x -\mu \|^2}{2z}}  \, dF_Z(z)\,\, dx\,.    
\end{eqnarray*}
Finally, part {\bf (c)} follows from pairing Theorem 2.1 of Strawderman (1974)  with part {\bf (b)} above, for the model $X'-\mu \sim SN_p(F_Z)$, with the finiteness conditions $E_0 \|X'X\|= p E(Z)< \infty$ and  $E_0 \|X'X\|^{-1} = \frac{1}{p-2} \, E(Z^{-1}) < \infty$, and with the upper bound on the constant $a$ given by $\frac{2}{E_0 \|X'X\|^{-1}}$.    \qed
\begin{remark}
\noindent  The dual loss in (\ref{reflectedsmnloss}) may be labelled as reflected scale mixture of normals analogously to reflected normal loss.  For the normal case with degenerate $W_1, W_2, V_1$ at $\sigma_Y^2$ and $\sigma_X^2$ respectively, 
part {\bf (a)} of Theorem \ref{smnl2loss} applied to the MRE estimator reduces to Corollary \ref{cplugin}, while part {\bf (b)} reduces to Corollary \ref{mreimprovements} with $Z$ in (\ref{Z}) degenerate at $\frac{\sigma_X^2 (\sigma_X^2 + 2 \sigma_Y^2)}{ 2(\sigma_X^2 +  \sigma_Y^2)}$.
\end{remark}
\noindent We pursue with some examples of applications of Theorem \ref{smnl2loss}.

\begin{example}  (Cases where both $H$ and $G$ are Gamma cdf's)
We illustrate some of the features of Theorem \ref{smnl2loss} for situations where
$W_1, W_2 \sim H \sim \hbox{Gamma}(\alpha_1,1)$, $V_1 \sim G \sim  \hbox{Gamma}(\alpha_2,1)$, with $\alpha_1 > p/2$ for $i=1,2$ which guarantees that $E(V_1^{-p/2}) < \infty$ and 
$E(W_1^{-p/2}) < \infty$ (see Definition \ref{definition1}).  We have $V_1+W_1 \sim F \sim \hbox{Gamma}(\alpha_1+\alpha_2,1)$ and $V_1+W_1+W_2 \sim J \sim \hbox{Gamma}(\alpha_1+2\alpha_2,1)$.  Part {\bf (d)} of Example \ref{bayesexamples} indicates that
\begin{equation*}
\hat{q}_{\hbox{mre}}(y;x) \,=\, q*p(y-x), \hbox{ with } \; q*p \sim SN_p(F)\,.
\end{equation*}
Theorem \ref{smnl2loss} tells us that $\hat{q}_{\hbox{mre}}(y;X)$ is inadmissible as an estimator of $q(y-\mu)$, $y \in \mathbb{R}^p$ for $p \geq 3$, and dominated by any 
$q*p(y-\hat{\mu}(X))$, where $\hat{\mu}(X')$ dominates $X'$ under squared error loss and for $X'-\mu \sim SN_p(F_Z)$ as given in (\ref{Z}).  The joint density of $(Z_1,Z_2)$ in
(\ref{Z}) becomes
\begin{equation*}
d\tau(z_1,z_2) \, = \, k \frac{z_1^{\alpha_1-1} \, z_2^{\alpha_1+2\alpha_2-2} \, e^{-(z_1+z_2)}}{(z_1+z_2)^{p/2}}  \; \mathbb{I}_{\mathbb{R}_+}(z_1) \, \mathbb{I}_{\mathbb{R}_+}(z_2)\,,
\end{equation*}
with $k = \frac{\Gamma(2\alpha_1+2\alpha_2-1)}{\Gamma(\alpha_1) \, \Gamma(\alpha_1+2\alpha_2-1) \, \Gamma(2\alpha_1+2\alpha_2-1-p/2)}$, and where we have used the identity
\begin{equation}
\nonumber
\int_{\mathbb{R}_+} \int_{\mathbb{R}_+}  \frac{z_1^{a-1} \, z_2^{b-1}}{(z_1+z_2)^{c}}
\, e^{-(z_1+z_2)} \, dz_1 dz_2 \, =\, \frac{\Gamma(a) \,\Gamma(b) \, \Gamma(a+b-c)}{\Gamma(a+b)}\,,
\end{equation}
for $a,b >0$, $a+b >c$.  Finally, using again the above identity, calculations yield
$$E_{\tau}(Z^{-1}) = E_{\tau}({Z_1}^{-1}) + E_{\tau}({Z_2}^{-1})\, = \, 
\frac{(\alpha_1+\alpha_2 -1) ( 2\alpha_1 + 2\alpha_2 -3)}{(\alpha_1+\alpha_2-1-p/4) \,
(\alpha_1-1) \, (\alpha_1+2\alpha_2-2)}\,,$$  
\end{example}
\noindent so that part {\bf (c)}'s subclass of dominating Baranchik predictive density estimators is explicitly determined with $0 <a \leq \frac{2 (p-2)}{E_{\tau}(Z^{-1})}$ and the above $E_{\tau}(Z^{-1})$.
 
\begin{example}  (Cases where the mixing distribution $G$ is lower bounded) 
Consider situations where either $X-\mu \sim N_p(0,\sigma_X^2 I_p)$ or, more generally, the scale parameter distribution for $X$ is bounded below by some known $a_X>0$ (i.e., $G^-(a_X)=0$ where $G^-$ is the left-hand limit of $G$).  With such an assumption, without any additional knowledge on $G$, one can obtain an upper bound for Theorem \ref{smnl2loss}'s $E(Z^{-1})$ and, hence, a lower bound for part {\bf (c)}'s upper limit 
$\frac{2(p-2)}{E(Z^{-1})}$.  Indeed, the lower bound assumption implies that $P(Z_1 \geq a_X)=1$, $P(Z_2 \geq a_X)=1$, $E(Z^{-1}) \leq \frac{2}{a_X}$, and $0<a<(p-2) a_X$ for Theorem \ref{smnl2loss}'s Baranchik-type estimators $\hat{\mu}_{a,r(\cdot)}$.   Similarly, if the mixing variance distribution $H$ for $Y$ is bounded below by some $a_Y>0$, the above bounds become  $E(Z^{-1}) \leq \frac{1}{a_X} + \frac{1}{a_X+2a_Y}$, and $0<a \leq \frac{(p-2) a_X (a_X+2a_Y)}{a_X+a_Y}$; with the degenerate case bringing us back to Example \ref{baranchikexample}. 
\end{example}

\subsection{Dominating estimators of a plug-in estimator by expanding the scale }
\label{scaleexpansionsection}
\noindent  We consider here scale mixtures of normals in (\ref{model}) as described in 
Definition \ref{definition1}:
\begin{equation}
\label{setupsmn}
X-\mu \sim SN_p(G)\,,\,\; Y-\mu \sim SN_p(H)\,.
\end{equation}
As in Section \ref{duality}, we focus on the class of estimators $\hat{q}_c(y;x)= \frac{1}{c^p} \, q(\frac{y-x}{c})$
of the density $q(y-\mu)$, $y \in \mathbb{R}^p$, with $c=1$ corresponding to the plug-in maximum likelihood estimator $q(y-x)$, and choices $c>1$ representing expansions of the known scale associated with the underlying density $q(y-\mu)$.  As shown in Theorem \ref{scalemixturesexpansion} and Remark \ref{remarkthatfollows}, we extend parts {\bf (b)} and {\bf (c)} of Corollary \ref{cor}, applicable to the normal case, which : {\bf (i)} provides the best estimator within the class of $\hat{q}_c$'s;  {\bf (ii)} shows the superiority of a subclass of estimators 
$\hat{q}_c$ over the plug-in $\hat{q}_1$ for a range $(1,c_1)$ (say) of values of $c$; {\bf (iii)} establishes that $\hat{q}_c$ dominates $\hat{q}_1$ universally for all $c>1$ for sufficiently large dimension $p$.   The next intermediate result follows from Lemma \ref{L2distance} and gives the loss incurred by using $\hat{q}_c$ to estimate $q(y-\mu), y \in \mathbb{R}^p$. 

\begin{corollary}
\label{scaleexpansionlemma}
Whenever finite, the $L_2$ distance between densities $\hat{q}_c(y;x)= \frac{1}{c^p} \, q(\frac{t-x}{c})$ and $q(t-\mu)$; $t,x,\mu \in \mathbb{R}^p$, $c>0$, is given by
\begin{equation}
\label{losshatqc}  (\frac{1}{c^p} +1) \, q*\bar{q}(0) \, - \, 2 q*\bar{h}(x-\mu)\, 
\end{equation}
where $h(t)= \frac{1}{c^p} \, q(\frac{t}{c})$ and $\bar{h}(t)=h(-t)$ for all $t$.  
\end{corollary}
\noindent {\bf Proof.}  With the evaluation $h*\bar{h}(0)= \frac{1}{c^p} \, 
q*\bar{q}(0)$, the result follows as Lemma \ref{L2distance} by setting $\mu_1=x, \mu_2=\mu$, and $f \equiv h\,$.                             \qed

\begin{theorem}
\label{scalemixturesexpansion}
Let $W_1,W_2 \sim H$, $V_1 \sim G$ be independently distributed.  Define 
$M_c=E[(V_1+W_1+c^2 \, W_2)^{-p/2}\,]$; $c>1$;  $N=E[(W_1+W_2)^{-p/2}\,]$.\footnote{The assumptions of Definition \ref{definition1} imply the finiteness of $M_c$ and $N$.}
\begin{enumerate}
\item[ {\bf (a)}]
The risk $R(\mu,\hat{q}_c)$ for estimating a normal scale mixture density $q(t-\mu)$, $t \in \mathbb{R}^p$, with $q \sim SN_p(H)$ under $L_2$ loss and for $X \sim SN_p(G)$, is constant in $\mu$ and given by
\begin{equation}
\label{smnriskqhatc}
R(\mu,\hat{q}_c)\,=\, \frac{1}{(2\pi)^{p/2}} \, \left( (1+\frac{1}{c^p})\,N\;-\; 2 M_c \right)\,;
\end{equation}
 
\item[ {\bf (b)}]  The optimal estimator among the class of estimators $\hat{q}_c$ is 
$\hat{q}_{c^*}$, where $c^*$ is the unique value of $c >1$ such that 
\begin{equation}
\label{cstar}
 2 c^{p+2} \, E[\frac{W_2}{(V_1+W_1+c^2 \, W_2)^{p/2+1}}]\, = \, N\,;
\end{equation}
\item[ {\bf (c)}] Estimators $\hat{q}_c$ with $1<c<c_1$ dominate the plug-in
 $\hat{q}_1$, where $c_1= \infty$  if $N \geq 2M_1$ and, otherwise,  $c_1$ is the unique value of $c >1$ such that $N(1-\frac{1}{c^p}) = 2 (M_1-M_c)\,.$
\end{enumerate}
\end{theorem}
\noindent {\bf Proof.}  {\bf (a)} The risk of $\hat{q}_c$ is given by 
\begin{equation}
\label{riska}
R(\mu,\hat{q}_c)\, =\, E_{\mu}\left[ (1+\frac{1}{c^p})\, q*q(0) \,-\, 2 \,q*h(X-\mu)\,  \right]\,,
\end{equation}
by taking the expected value of (\ref{losshatqc}) with respect to $X-\mu \sim SN_p(G)$, and since $\bar{q} \equiv q$ and $\bar{h} \equiv h$ by spherical symmetry of $q$ and thus also $h$. 
By making use of Lemma \ref{convolutionscalemixture}, we have for the scale mixture of normals density, $q  \sim SN_p(H)$, the convolutions $q*q \sim SN_p(F_1)$ and $q*h \sim SN_p(F_c)$, $F_c$ being the cdf of $W_1+c^2W_2$ for any $c^2$. With the
above, we obtain 
\begin{equation}
\label{qq0}  q*q(0) \,=\, (2\pi)^{-p/2}\, N\,.
\end{equation}
Furthermore, we have 
\begin{eqnarray}
\nonumber  E_{\mu} [\,q*h(X-\mu)\,] &=& \int_{\mathbb{R}^p} \int_{\mathbb{R}_+} (2\pi v)^{p/2} 
\,\left(\int_{\mathbb{R}_+}(2\pi t)^{p/2} \, e^{- \frac{\|x-\mu\|^2}{2t}} \, dF_c(t) \right) \,  e^{- \frac{\|x-\mu\|^2}{2v}} \, dG(v) \, dx \\
\nonumber \, &=& \int_{\mathbb{R}_+} \int_{\mathbb{R}_+} (4\pi^2 vt)^{-p/2}\, 
\left(\int_{\mathbb{R}^p} e^{- \frac{\|x\|^2}{2} (\frac{1}{t} + \frac{1}{v}}) \,  dx \right) \, dF_c(t) dG(v) \\
\nonumber\, &=&   \int_{\mathbb{R}_+} \int_{\mathbb{R}_+} (4\pi^2 vt)^{-p/2}\, 
 (\frac{2\pi\, vt}{v+t})^{p/2} \; dF_c(t) dG(v) \\
\label{q*h} \,&=& (2\pi)^{-p/2} \, E[(V_1+W_1+c^2 \, W_2)^{-p/2}]\,=\, (2\pi)^{-p/2}\, M_c\,.
\end{eqnarray}
Finally, the given expression for $R(\mu,\hat{q}_c)$ in (\ref{smnriskqhatc}) follows from
(\ref{riska}), (\ref{qq0}), and (\ref{q*h}). \\

\noindent {\bf (b)}  It is easy to see from (\ref{smnriskqhatc}) that
\begin{equation*}
\frac{(2\pi)^{p/2} \, c^{p+1}}{p} \; \frac{\partial}{\partial c} R(\mu,\hat{q}_c)\,=\,
-N + 2 \, E[\,W_2 (\frac{c^2}{V_1+W_1+c^2 \, W_2})^{p/2+1}\,]\,.  
\end{equation*}
Evaluated at $c=1$, the above is negative, while it is positive evaluated at $c \to \infty$.  Moreover, since the above is increasing as a function of $c \in [1,\infty)$, we have that $\frac{\partial}{\partial c} R(\mu, \hat{q}_c)$ changes signs once from $-$ to $+$ as $c$ increases on $[1,\infty)$.     \\
 
\noindent {\bf (c)} Given that, as a function of $c$, $R(\mu, \hat{q}_c)$ is strictly decreasing for $1 \leq c < c^*$, and strictly increasing for $c>c^*$, we have indeed $R(\mu,\hat{q}_c) < R(\mu,\hat{q}_1)$ for all $c>1$ as soon
as $\lim_{c \to \infty}\, R(\mu,\hat{q}_c) \leq  R(\mu,\hat{q}_1) \Longleftrightarrow 2M_1\leq N $.  Otherwise, we have $R(\mu,\hat{q}_c) < R(\mu,\hat{q}_1)$ for $c<c_1$ with $R(\mu,\hat{q}_{c_1}) =  R(\mu,\hat{q}_1) \Longleftrightarrow N(1-\frac{1}{c^p}) = 2 (M_1-M_c)\,.$    \qed \\

\noindent  
\begin{remark}
\label{remarkthatfollows}
The above Theorem is presented for fixed $p$, but there also implications for varying $p$ analogously to part {\bf (c)} of Corollary \ref{cor} established for normal models.
Indeed, assuming all the inverse moments associated with $H$ and $G$ exist, which  guarantees the finiteness of $N$ and $M_1$ for all $p \geq 1$, it is inevitable that the interval of values of $c$ such that $\hat{q}_c$ dominates $\hat{q}_1$ is given by $(1,\infty)$ for large enough $p \geq p_0$. This is justified by the fact that if $N \geq 2 M_1$ for a given $p_0$ (which can be shown to exist), i.e.,  
$$  E (\frac{1}{(W_1+W_2)^{p_0/2}}) \geq 2  E (\frac{1}{(V_1+ W_1+W_2)^{p_0/2}}\,)\,,$$
then we must also have
$$  E (\frac{1}{(W_1+W_2)^{(1+p_0)/2}}) \geq 2  E (\frac{1}{(V_1+ W_1+W_2)^{(1+p_0)/2}}\,)\,,$$
telling us that $N \geq 2 M_1$ for all $p >p_0$.
\end{remark}
 
\begin{example}
Theorem \ref{scalemixturesexpansion} and Remark \ref{remarkthatfollows} apply in the normal case with $P(W_i=\sigma_Y^2)=1$ for $i=1,2$ and $P(V_1=\sigma_X^2)=1$, and it is readily verified that results in parts {\bf (b)} and 
{\bf (c)} of Corollary \ref{cor} follow.  \\
As a further illustration, consider situations where $W_1, W_2, V_1$ share the same distribution with $P(a_1 \leq W_1 \leq a_2)=1$ for some $0 <a_1 \leq a_2 < \infty$.
By setting
$T=V_1+W_1 =^d W_1+W_2$ and $D=W_2/T$, equation (\ref{cstar}) may be expressed as
\begin{equation*}
 c^{p+2}  \, E^T \left( T^{-\frac{p}{2}} E^{D|T} (\frac{2D}{(1+c^2D)^{p/2+1}} )\right)   =  E^T \left( T^{-\frac{p}{2}} \right)\,,
\end{equation*}
so that the condition 
\begin{equation}
\label{c(t)}
C(t) \, = \,  c_0^{p+2}\, E^{D|T=t} (\frac{2D}{(1+c_0^2D)^{p/2+1}} )  \leq 1\,,
\end{equation}
for all $t \in [2r_1, 2r_2]$, and for some $c_0$, implies $c^* \geq c_0$.  Using the covariance inequality $\hbox{Cov}(f_1(D), f_2(D)) \leq 0$ for increasing $f_1$ and decreasing $f_2$, as well as the property $E(D|T)=1/2$ which is a consequence of the iid assumption on $W_1, W_2, V_1$, we obtain 
\begin{equation*}
C(t) \leq (c_0)^{p+2} \, E^{D|T} (\frac{1}{(1+c_0^2D)^{p/2+1}} )\,.
\end{equation*}
Now, with the bounded support assumption setting $\beta = \frac{r_1}{r_1+r_2}$, we have $P(D \geq \beta|T=t) =1$ and $C(t) \leq (\frac{c_0^2}{1+c_0^2 \beta})^{p/2+1}$ for all $t$.  Finally, setting this upper bound on $C(t)$ equal to $1$, we obtain the lower bound
$$ {c^*}^2 \geq 1 + \frac{r_1}{r_2}\,.$$ \\

\noindent  Turning now to a value $p_0$ such that all $\hat{q}_c$ with $c>1$ dominate the plug-in $\hat{q}_1$ for all $p \geq p_0$, Theorem \ref{scalemixturesexpansion}'s condition $N \geq 2M_1$ may be written as
\begin{equation*}
N \geq 2M_1  \Longleftrightarrow  E(T^{-\frac{p}{2}}) \geq 2 E^T \left(T^{-\frac{p}{2}} \;  E^{D|T} (1+D)^{-\frac{p}{2}}  \right)\,,
\end{equation*}
which becomes satisfied as soon as $E^{D|T=t} (1+D)^{-\frac{p}{2}} \leq \frac{1}{2}$ for all $t$.  Finally, with $P(D \geq \beta|T=t) =1$, we conclude that all $\hat{q}_c$ with $c>1$ dominate the plug-in $\hat{q}_1$ for all $p \geq p_0= \frac{\log 4}{\log (1+\beta)}\,. $
\end{example}
 
\section{Integrated $L_1$ loss and plug-in estimators}
\label{L1section}
The results of this section apply, as in Section \ref{extension}, to the spherically symmetric set-up in (\ref{ss}), but relate to integrated $L_1$ loss.  We focus on the performance of plug-in estimators 
$q_Y(\|y- \hat{\mu}(X) \|^2)$, $y \in \mathbb{R}^p$, with $\hat{\mu}(X)$ an estimator of $\mu$.  We capitalize 
on an explicit representation for the $L_1$ distance (Lemma \ref{ssdistance}) between two densities of the same spherically symmetric 
family to establish that our predictive density estimation problem for plug-in estimators is dual to a point estimation problem for the same plug-in estimators
under a loss which is a concave function of $\| \hat{\mu} -\mu\|^2$ (Corollary \ref{pluginL1}).  
As in Sections \ref{normalsection} and \ref{extension}, using Stein estimation results and techniques applicable to such concave losses, we establish the inadmissibility of plug-in densities $q_Y(\|y- X \|^2)$ for $p \geq 4$ and obtain dominating predictive density estimators.  
In subsection \ref{L1smnsection}, we provide further specific developments for scale mixtures of normals $p_X$ and $q_Y$, which includes of course the normal case.
Finally, in subsection \ref{L1restrictionssection}, for univariate situations where $\mu$ is either restricted to an interval $(a,b)$ or restricted to a left-bounded interval $(a,\infty)$, we proceed as in Corollary \ref{ms2005} to show that the plug-in density estimator $q_Y(|y- \hat{\mu}_{\pi_U}(X) |^2)$ dominates the plug-in $q_Y(|y- X |^2)$ for log-concave $p_X(x^2)$ in $x$, where $\hat{\mu}_{\pi_U}(X)$ is the Bayes point estimator of $\mu$ associated with a uniform prior on the restricted parameter space and the given dual loss.  

\subsection{An identity for $L_1$ distance and general dominance results of plug-in predictive density estimators}
\label{L1identitysection}
\noindent We begin with a useful $L_1$ distance identity which is also of independent interest.  

\begin{lemma}
\label{ssdistance}
Let $Y=(Y_1, \ldots, Y_p)'$ be a spherically symmetric distributed random vector with unimodal, Lebesgue density $q_Y(\|y-\mu\|^2)$; $y \in \mathbb{R}^p$.
Then for any $\mu_1, \mu_2 \in \mathbb{R}^p$, the $L_1$ distance between $f_{\mu_1}$ and $f_{\mu_2}$ is given by
\begin{equation}
\rho_{L_1}= \int_{\mathbb{R}^p} |q_Y(\|y-\mu_1\|^2) - q_Y(\|y-\mu_2\|^2)| \,dy \,=\, 4 F(\frac{\|\mu_1-\mu_2\|}{2}) -2\,,
\end{equation} 
where $F(t)=P_0(Y_1 \leq t)$, $t \in \mathbb{R}$, is the cumulative distribution function of $Y_1$ when $\mu_1=0$.  
\end{lemma}

\begin{remark}
\label{concave}
This result was given by Das Gupta and Lahiri (2012) for the normal case.  An existing reference for the general case seems likely to us, but we could not find such a reference.  Observe that the distance $\rho_{L_1}$ is always a concave function of $\|\mu_1-\mu_2\|$ on $(0,\infty)$ since $F'$ is unimodal, and also of $\|\mu_1-\mu_2\|^2$ given that $F$ is increasing. 
\end{remark}
\noindent {\bf Proof.}
 We have $q_Y(\|y-\mu_1\|^2) \geq q_Y(\|y-\mu_2\|^2) \Leftrightarrow \|y-\mu_1\|^2 \leq \|y-\mu_2\|^2 \Leftrightarrow L(y) \leq 0$, where $L(y)= (\mu_2-\mu_1)'y +
\frac{\|\mu_1\|^2 - \|\mu_2\|^2}{2}$.  Setting $A=\{y \in \mathbb{R}^p: L(y) \leq 0  \}$, we obtain splitting the integration on $A$ and its complement $A^c$
\begin{eqnarray}
\nonumber \rho_{L_1} &=& P_{\mu_1}(Y \in A) - P_{\mu_2}(Y \in A) + P_{\mu_2}(Y \in A^c) -  P_{\mu_1}(Y \in A^c) \\
\label{rho} \,& = & 2 \, \{P_{\mu_1}(L(Y) \leq 0) + P_{\mu_2}(L(Y) \geq 0) -1   \}\,.
\end{eqnarray}
Observe that $L(Y)$ is a linear function of the spherically symmetric distributed $Y$.  For such linear functions, we have (e.g., Muirhead, 2005) 
$$  \frac{(l'Y+k) - (l'\mu+k)}{\|l\|}  \sim F\,,$$
for all $l \in \mathbb{R}^p - \{0\}, k \in \mathbb{R}^1$.  We thus obtain $P_{\mu_1}(L(Y) \leq 0) = F\left(-(L(\mu_1))   \right) = 
F(\frac{\|\mu_1 -\mu_2\|}{2}).$  Similarly, we obtain $P_{\mu_2}(L(Y) \geq 0) = F(\frac{\|\mu_1 -\mu_2\|}{2})$, and the desired expression for $\rho_{L_1}$
follows from (\ref{rho}).  \qed

\begin{corollary}
\label{pluginL1}
For estimating an unimodal spherically symmetric Lebesgue density $q_Y(\|y-\mu\|^2)$, $t \in \mathbb{R}^p$, under integrated $L_1$ loss and based on
$X \sim p_X(\|x-\mu\|^2)$, the frequentist risk of the plug-in density estimator $q_Y(\|y-\hat{\mu}(X)\|^2)$ is equal to the frequentist risk of the point estimator $\hat{\mu}(X)$ of $\mu$ under loss $4 F(\frac{\|\hat{\mu}-\mu\|}{2}) -2\,$, with $F$ being the common marginal cdf associated with $q_Y$. Consequently, $q_Y(\|y-\hat{\mu}_1(X)\|^2)$ dominates $q_Y(\|y-\hat{\mu}_2(X)\|^2)$ iff $\hat{\mu}_1(X)$ dominates $\hat{\mu}_2(X)$ under loss $2 F(\frac{\|\hat{\mu}-\mu\|}{2}) -1\,$.
\end{corollary}
\noindent {\bf Proof.}  This is a direct consequence of Lemma \ref{ssdistance}.  \qed \\

\noindent   Since the dual problem described above involves loss functions $l(\|d-\mu\|^2)$  with $l(t)=2F(\frac{\sqrt{t}}{2}) -1 $ being concave (see Remark \ref{concave}), we consider using Stein estimation techniques and results for such concave losses ( Brandwein and Strawderman 1991, 1981); Brandwein, Ralescu and Strawderman, 1993), along with  Corollary \ref{pluginL1}, to obtain dominating estimators of the plug-in density $q_Y(\|y-X\|^2)$, $y \in \mathbb{R}^p$, which we now proceed to do, elaborate on, and illustrate.  For what follows, we denote $f$ as the density of $\|X-\mu\|$ under $p_X$ and we recall that $f(t)= \frac{2\pi^{p/2}}{\Gamma(p/2)} \, t^{p-1} \, p_x(t^2)\,$ (e.g., Muirhead, 2005).  Here is an adaptation of Theorem 2.1 of Branwein and Strawderman (1981) applicable to Baranchik type estimators, and followed by related inferences for improving on plug-in density estimators under $L_1$ loss.  

\begin{theorem} (Brandwein and Strawderman, 1981)
\label{bs1981}
Let $X$ have a spherically symmetric distribution with density $p_X(\|x-\mu\|^2)$, $x\in \mathbb{R}^p $, with respect to $\sigma$-finite measure $\nu$.  For $p \geq 4$ and for estimating $\mu \in \mathbb{R}^p$ under loss $l(\|d-\mu\|^2)$ with $l$ non-decreasing and concave on $\mathbb{R}_+$, estimators $\hat{\mu}_{a,r(\cdot)}(X)=(1 - a \frac{r(X'X)}{X'X}) X$ dominate $X$, and are thus minimax, provided:
\begin{enumerate}
\item[ (i)] $0 \leq r(\cdot) \leq 1$ and $r(\cdot) \neq 0$;
\item[ (ii)] $r(t)$ is non-decreasing for $t>0$;
\item[ (iii)] $r(t)/t$ is non-increasing for $t>0$; 
\item[ (iv)]  $0< E_{p_X} l'(\|X-\mu\|^2) < \infty$;
\item[ (v)]  $0 < a \leq \frac{2p}{p-2} \, \frac{1}{E_h (R^{-2})}\,$,
where $h(s)$ is a density on $\mathbb{R}_+$  proportional to $l'(s^2) \, f(s) = \frac{2\pi^{p/2}}{\Gamma(p/2)} \, l'(s^2)  s^{p-1} \, p_x(s^2) \,$. 
\end{enumerate}
\end{theorem}


\noindent This now follows from Corollary \ref{pluginL1} and Theorem \ref{bs1981}.  \qed

\begin{corollary}
\label{implication}
For estimating a unimodal spherically symmetric Lebesgue density $q_Y(\|y-\mu\|^2)$, $y,\mu \in \mathbb{R}^p$ and $p \geq 4$, under integrated $L_1$ loss and based on
$X \sim p_X(\|x-\mu\|^2)$, a plug-in Baranchik density estimator $q_Y(\|y-\hat{\mu}_{a,r(\cdot)}(X)\|^2)$  with $\hat{\mu}_{a,r(\cdot)}(X)=(1 - a \frac{r(X'X)}{X'X}) X$
dominates the plug-in $q_Y(\|y-X\|^2)$ provided conditions (i), (ii), and (iii) of Theorem \ref{bs1981} are satisfied as well as:
\begin{enumerate}
\item[ (iv)']  $0 < E_{p_X} (\frac{q_Y(\|X-\mu\|^2/16)}{\|X-\mu\|})< \infty$;
\item[ (v)']  $0 < a \leq \frac{2p}{p-2} \, \frac{\int_{(0,\infty)} u^{\frac{p-3}{2}} \, p_X(u) \, F'(\frac{u}{4})\, d\nu(u) }
{\int_{(0,\infty)} u^{\frac{p-5}{2}} \, p_X(u) \,F'(\frac{u}{4})\, d\nu(u)}\,$.
\end{enumerate}
\end{corollary}
\noindent  {\bf Proof.}  This follows from Corollary \ref{pluginL1} and Theorem \ref{bs1981} with $l(u)=2F(\frac{\sqrt{u}}{2}) -1$ and $l'(u)=\frac{F'(\frac{\sqrt{u}}{2})}{2\sqrt{u}}$.

\begin{remark}
In our setup, the model density $q_Y$ determines the loss $l$ via Lemma \ref{ssdistance} and is thus taken to be unimodal and Lebesgue.  On the other hand, there no restrictions on $p_X$ other than risk-finiteness for the estimators $\hat{\mu}_{a,r(\cdot)}(X)$.  Condition (iv') is weak.  For instance, it is satisfied when both the densities $q_Y$ and $p_X$ are bounded.   The upper bound for the multiplier $a$ in the estimator $\hat{\mu}_{a,r(\cdot)}(X)$ in condition (v') depends on both $q_Y$ and $p_X$.   
\end{remark}  
\noindent Here is an evaluation for the particular case when both $p_X$ and $q_Y$ are normal densities.
\begin{example} (normal case)
\label{normalL1}  For the normal case (\ref{normalmodel}) with $q_y(u) = (2\pi \sigma_Y^2)^{-p/2} \, e^{-u/2\sigma_Y^2}$ and $p_X(u) = (2\pi \sigma_X^2)^{-p/2} \, e^{-u/2\sigma_X^2}$, Corollary \ref{implication} applies with (iv') satisfied and (v') specializing to 
\begin{equation}
\label{cutoff}
0 < a \leq \frac{2p}{p-2} \,  \frac{\int_{0,\infty} u^{\frac{p-3}{2}} \,e^{-u/2\sigma_X^2} e^{-u/8\sigma_Y^2}  \, du }
{\int_{0,\infty} u^{\frac{p-5}{2}} \, \,e^{-u/2\sigma_X^2} e^{-u/8\sigma_Y^2}  \, du}\,= \, \frac{(p-2) (p-3)}{p} \, \frac{8 \sigma_X^2 \sigma_Y^2}{\sigma_X^2+4\sigma_Y^2}\,. 
\end{equation}
 We point out that a simultaneous dominance result is available for a family of $p_X$ models by taking the infimum with respect to $p_X$ on the rhs of (v').  For the normal case, if we have for instance $X \sim N_p(\mu, \sigma_X^2 I_p)$ with $\sigma_X^2$ unknown, but known to bounded below by $a_X>0$, then simultaneous dominance occurs for all such $p_X$'s by taking $0 < a \leq \frac{(p-2) (p-3)}{p} \, \frac{8 a_X \sigma_Y^2}{a_X+4\sigma_Y^2}\,$.
\end{example}

\subsection{Improvements for scale mixture of normals}
\label{L1smnsection}
Further developments for scale mixtures of normals are provided in this section and lead to wider classes of dominating estimators than those given by Corollary \ref{implication}.  We revisit this latter corollary for situations in (\ref{model}) where
\begin{equation}
\label{scalemixturerepresentation}
X - \mu \sim  SN_p(G)\,  \; Y - \mu \sim  SN_p(H).
\end{equation}
We define $Z$ as a random variable, $F_Z$ as its cdf, and $\tau$ as a bivariate cdf such that 
\begin{equation}
\label{Z}
Z =^d \frac{4Z_1Z_2}{Z_1+4Z_2}\,, \hbox{ with } (Z_1,Z_2) \sim d\tau(z_1,z_2) \propto \frac{z_2^{p/2-1}}{(z_1+4z_2)^{p/2}} \, dG(z_1) \, dH(z_2) \,.
\end{equation}

\begin{theorem}
\label{L1smn}  Let $X \sim p_X(\|x-\mu\|^2) $ and $Y \sim q_Y(\|y-\mu\|^2)$, $x,y,\mu \in \mathbb{R}^p$,  be scale mixtures of normals as in (\ref{scalemixturerepresentation}) and consider estimating $q_Y(\|y-\mu\|^2)$ under integrated $L_1$ loss based on $X$. 

\begin{enumerate} 
\item[ {\bf (a)}] For $p>1$,\footnote{For $p=1$, the density in (\ref{p*}) is not well defined.} the plug-in density estimator $q_Y(\|y-\hat{\mu}(X)\|^2)$  dominates the plug-in density $q_Y(\|y-X\|^2)$ provided $\hat{\mu}(X')$ dominates
$X'$ under loss $\| \hat{\mu} - \mu \|^2$ and for $X' \sim p^*(\|x-\mu\|^2)$, with
\begin{equation}
\label{p*}
p^*(\|s\|^2) = \frac{K}{\|s\|}  \, \int_{(0,\infty)}  (2\pi z)^{-p/2} \, 
e^{-\frac{\|s\|^2}{2z}} \, dF_Z(z)\,, s \in \mathbb{R}^p\,.
\end{equation}
\item[ {\bf (b)}]
In particular, a plug-in Baranchik density estimator $q_Y(\|y-\hat{\mu}_{a,r(\cdot)}(X)\|^2)$,  with $\hat{\mu}_{a,r(\cdot)}(X)=(1 - a \frac{r(X'X)}{X'X}) X$, dominates the plug-in density $q_Y(\|y-X\|^2)$ provided conditions (i), (ii), and (iii) of Theorem \ref{bs1981} are satisfied, $p \geq 4$,  the expectations $E(Z^{1/2})$ and $E(Z^{-3/2})$ are finite, and  $0<a \leq 2 (p-3) \frac{E(Z^{-1/2})}{E(Z^{-3/2})}\,.$
\end{enumerate}
\end{theorem}
\noindent{\bf Proof.}  {\bf (a)} We apply Corollary \ref{pluginL1}.  We thus seek conditions for which the difference in risks $\Delta(\hat{\mu},\mu) =
E_{\mu} \left[\, l (\|\hat{\mu}(X) - \mu \|^2) \, - \, l (\|X - \mu \|^2)   \,\right] $ is less than $0$, where $l(\|\hat{\mu}-\mu \|^2) = F(\frac{\|\hat{\mu}-\mu\|}{2}) -1\,$.  We apply the inequality $l(s) - l(t) < l'(t) (s-t)$ for strictly concave $l$ and $s \neq t$, which implies for the difference in losses that
\begin{equation}
\label{differencelosses}
l(\|\hat{\mu}(x)-\mu \|^2) \, - \,  l (\|x - \mu \|^2) <
l'(\|x - \mu \|^2) \, \left(\|\hat{\mu}(x)-\mu \|^2 - \|x-\mu \|^2  \right)\,,
\end{equation}
for all $x, \mu \in \mathbb{R}^p$ such that $x \neq \hat{\mu}(x)$.  Observe that
\begin{eqnarray}
\nonumber  l'(\|x - \mu \|^2) &=&  \frac{1}{2 \|x -\mu \|} \, F'(\frac{\|x -\mu \|}{2}) \\
\label{lprime}\, & =& \frac{1}{2 \|x -\mu \|} \,   \int_0^{\infty}  (2\pi w)^{-1/2} \, 
e^{-\frac{\|x -\mu \|^2}{8w}} \, dH(w)\,,
\end{eqnarray}
since the marginal distributions associated with a scale mixture of normals as in (\ref{scalemixturerepresentation}) are themselves univariate scale mixtures of normals with the same mixing distribution.  \footnote{It is not the case that spherically symmetric distributions share a similar consistency property, but it is true for scale mixtures of normals distributions (e.g., Kano, 1994) .}
Now, using (\ref{differencelosses}) and (\ref{lprime}), it follows that
\begin{eqnarray*}
\Delta (\hat{\mu}, \mu) &<&  E_{\mu}^X  \left[
\frac{ \left(\|\hat{\mu}(X)-\mu \|^2 - \|X-\mu \|^2  \right)}{2 \|X -\mu \|} \,  \int_0^{\infty} (2\pi w)^{-1/2} \, 
e^{-\frac{\|X -\mu \|^2}{8w}} \, dH(w) \right] \\
\, &=& \int_{\mathbb{R}^p} \frac{ \left(\|\hat{\mu}(x)-\mu \|^2 - \|x-\mu \|^2  \right) }{2 \|x -\mu \|} \,\int_0^{\infty} \int_0^{\infty} 
\frac{(2\pi)^{-\frac{p+1}{2}}}{(wv^p)^{\frac{1}{2}}}
 \, e^{-\frac{\|x -\mu \|^2}{(\frac{8wv}{v+4w})}}  dG(v) dH(w) \, dx \\
\, &\propto & \int_{\mathbb{R}^p} \frac{ \left(\|\hat{\mu}(x)-\mu \|^2 - \|x-\mu \|^2  \right) }{2 \|x -\mu \|} \,\int_0^{\infty} \int_0^{\infty}
(\frac{8\pi wv}{v+4w})^{-\frac{p}{2}} \, e^{-\frac{\|x -\mu \|^2}{(\frac{8wv}{v+4w})}} \,d\tau(v,w) \, dx \\
&\propto& \int_{\mathbb{R}^p} \frac{ \left(\|\hat{\mu}(x)-\mu \|^2 - \|x-\mu \|^2  \right) }{2 \|x -\mu \|} \,  \int_0^{\infty} (2\pi z)^{-p/2} \, e^{-\frac{\|x -\mu \|^2}{2z}} \, dF_Z(z)\,, dx\,,
\end{eqnarray*}
which establishes part {\bf (a)}.

\noindent  {\bf (b)}  We show below that the density $p^*(\|s\|^2)$ is a scale mixture of normals.  This permits us to apply the dominance result of Strawderman (1974) for Baranchik estimators satisfying conditions (i), (ii), (iii) of Theorem \ref{bs1981}, in cases where both $E_0^{p^*}\|X\|^2$ and $E_0^{p^*}\|X\|^{-2}$ are finite, and for $0<a \leq 2/(E_0^{p^*}\|X\|^{-2})$.  The finiteness conditions are satisfied for $p \geq 4$ and with the finiteness of $E(Z^{1/2})$ and $E(Z^{-3/2})$, and a calculation yields
\begin{eqnarray*}
E_0^{p^*}\|X\|^{-2} &=& \frac{\int_{\mathbb{R}^p} \frac{1}{\|x\|^{3/2}}
\int_0^{\infty} (2\pi z)^{-p/2} \, e^{-\frac{\|x\|^2}{2z}} \, d\tau(z) \, dx}{\int_{\mathbb{R}^p} \frac{1}{\|x\|^{1/2}}
\int_0^{\infty} (2\pi z)^{-p/2} \, e^{-\frac{\|x\|^2}{2z}} \, d\tau(z) \, dx } \\
\, &=&  \frac{\int_0^{\infty} \int_{\mathbb{R}^p} \frac{1}{\|x\|^{3/2}}
 (2\pi z)^{-p/2} \, e^{-\frac{\|x\|^2}{2z}} \, dx \, d\tau(z)}{\int_0^{\infty} \int_{\mathbb{R}^p} \frac{1}{\|x\|^{1/2}}
(2\pi z)^{-p/2} \, e^{-\frac{\|x\|^2}{2z}} \,dx \, d\tau(z)} \\
\, &=&  \frac{1}{(p-3)} \frac{E(Z^{-3/2})}{ E(Z^{-1/2})}\,,
\end{eqnarray*}
using expectation expressions for a central $\chi_p^2$ distribution.  This yields the desired result.  It remains to show that $p^*(\|s\|^2), s \in \mathbb{R}^p$, is a scale mixture of normals density.  Recall that, in general, a spherically symmetric density $f(\|t-\mu\|^2)$ is a scale mixture of normals if and only if $f$ is completely monotone on $(0,\infty)$, i.e., $(-1)^k \, 
f^{(k)}(t) \geq 0$ for $t>0$ and $k=0,1,2, \ldots$ (e.g., Berger, 1975).  Since 
both $t^{-1/2}$ and $\int_{(0,\infty)}  (2\pi z)^{-p/2} \, 
e^{-\frac{t}{2z}} \, d\tau(z)\,$ are completely monotone, it follows that their product is completely monotone (e.g., Feller, 1966, page 417) and that the density in (\ref{p*}) is indeed a scale mixture of normals.  \qed

\begin{example} (normal case)
In the normal case (\ref{normalmodel}) which arises as a particular case of (\ref{scalemixturerepresentation})
for degenerate $V,W$, we obtain that $Z$ in (\ref{Z}) is also degenerate with
$P(Z= z_0)=1$, with $z_0=\frac{4\sigma_X^2 \sigma_Y^2}{\sigma_X^2 + 4 \sigma_Y^2}$.  In this case as well, we obtain
$$  p^*(\|s\|^2) \, \propto \frac{1}{\|s\|}  (2\pi z_0)^{-p/2} \, e^{- \frac{-\|s\|^2}{2z_0}} \,,\;\;\; $$
which is the density of a Kotz distribution (see for instance Nadarajah, 2003).  By virtue of part (a) of Theorem \ref{L1smn}, minimax or dominance results applicable to this particular Kotz distribution generate plug-in $N_p(\hat{\mu}(X), \sigma_Y^2 I_p)$ density  estimators (such as those in part b) which dominate the
plug-in density of a $N_p(X, \sigma_Y^2 I_p)$ under $L_1$ loss. 
The cut-off point in part (b) reduces to $ 2(p-3) z_0 =  \, \frac{8 (p-3) \sigma_X^2 \sigma_Y^2}{\sigma_X^2+4\sigma_Y^2}\,$.  In comparison to Corollary \ref{pluginL1}'s cutoff point given in (\ref{cutoff}), the cut-off point here is larger by a multiple of $p/(p-2)$. 

\end{example}

\subsection{Improvements in the case of univariate parametric restrictions}
\label{L1restrictionssection}
We briefly expand on dominance results applicable to univariate ($p=1$) cases where $\mu$ is either restricted to an interval $(a,b)$ or a left-bounded interval $(a,\infty)$.  Combining  Corollary \ref{pluginL1}'s duality with point estimation loss $2 F(\frac{|\hat{\mu}-\mu|}{2}) -1$, which is a strictly bowled shaped function of $|\hat{\mu}-\mu|$ on $\mathbb{R}$, with findings of Marchand and Strawderman (2005), we derive an $L_1$ analog of Corollary \ref{ms2005} for estimating an 
univariate density $q_Y(|y-\mu|^2)$ based on $X \sim p_X(|x-\mu|^2)$ for cases (such as the normal case) where the family of densities for $X$ has an increasing monotone likelihood ratio (or equivalently $p_X(t^2)$ is log concave in $t \in \mathbb{R}_+$).   

\begin{corollary}
\label{ms2005L1}
For estimating an unimodal and univariate symmetric Lebesgue density $q_Y(|y-\mu|^2)$, $y \in \mathbb{R}$, $\mu \in (a,b)$ (or $\mu \in (a,\infty)$) under integrated $L_1$ loss and based on $X \sim p_X(|x-\mu|^2)$ with $p_X(t^2)$ logconcave, the plug-in density estimator $q_Y(|y-\hat{\mu}_U(X)|^2)$ with $\hat{\mu}_U(X)$ the Bayes estimator of $\mu$ with respect to the uniform prior on $(a,b)$ (or on $(a,\infty)$) dominates the plug-in density estimator  $q_Y(|y-X|^2)$.
\end{corollary}
\noindent  {\bf Proof.}  The result follows directly from Corollary \ref{pluginL1} and results in Marchand and Strawderman (2005).  \qed

\section{Concluding Remarks}
We have investigated the frequentist risk performance of various predictive density estimators under both integrated $L_2$ and $L_1$ loss.  For multivariate normal models, we have established a connection between the $L_2$ risk of plug-in type estimators and point estimation risk under reflected normal loss.  Paired with Stein estimation techniques and results for estimating a multivariate normal mean under loss which is a concave function of
the squared error $\|\hat{\mu}-\mu\|^2$, we establish the inadmissibility of the minimum risk equivariant (MRE) density estimator and obtain dominating predictive density estimators for three dimensions or more.  The duality is further exploited to obtain improvements of the benchmark MRE density estimator in the presence of restrictions on the underlying mean parameter.  
\noindent We have also analyzed the performance of scale expansion plug-in density estimators $N_p(\hat{\mu}(X), c^2 \sigma_Y^2 I_p)$ with varying $c^2$, obtaining notably instances (i.e., large enough dimension $p$ and $\hat{\mu}(x)=ax$ with $0<a \leq 1$) where all scale expansions $c^2>1$ improve uniformly on $c^2=1$. \\

\noindent  For scale mixtures of multivariate normal observables, we have obtained analogous 
developments with regards to the MRE density estimator by making use of a general $L_2$ distance identity, including its inadmissibility and the determination of explicit improvements, in general for three of more dimensions.  As well, we obtain improvements on the plug-in maximum likelihood estimator by scale expansion.   Finally, we have considered $L_1$ integrated loss and spherically symmetric observables, and, via an $L_1$ distance identity, obtained dominating estimators of a benchmark plug-in density estimator, in general for four dimensions or more.  \\

\noindent  In summary, the findings of this paper provide fundamental identities and results for assessing the efficiency of predictive density estimators of multivariate observables for both $L_2$ and $L_1$ integrated losses.  The main themes, simplified somewhat, revolve about the inefficiency of MRE estimators in high enough dimensions and about the inefficiency of plug-in estimators by either improving on the plug-in for a dual point estimation loss or expanding the scale.  Developments for such models with unknown scale represents one of several challenging and interesting problems worthwhile pursuing.   

\section{Appendix}
\label{Appendix}
{\bf A.1.  Proof of the minimaxity in Proposition \ref{mreprop}.}

\noindent We proceed as in Girshick and Savage (1951).
For $\mu=(\mu_1, \ldots, \mu_p)'$, let $A_k=\{ \mu | \ | \mu_i| < k/2, i=1, \ldots, p \}$ for $k=1, 2,\ldots$, and 
consider the sequence of prior distributions given by 
$$
\pi_k(\mu) =
\left\{
\begin{array}{cl}
k^{-p} & {\rm if} \ \mu\in A_k\\
0 & {\rm otherwise},
\end{array}
\right.
$$
which yields the Bayes estimators
$$
{\hat q}_k^\pi(y|x) = \int_{A_k} q(y-a) \, p(x-a) \,da / \int_{A_k} p(x-a) \, da
$$
with the Bayes risk function
$$
r_k(\pi_k, {\hat q}_k^\pi) =
{1\over k^p}\int_{A_k} \int\int \{ q(y-\mu)-{\hat q}_k^\pi(y|x)\}^2 \, dy \, p(x-\mu) \,dx \, d\mu.
$$
Since $r_k(\pi_k,{\hat q}_k^\pi)\leq r_k(\pi_k,{\hat q}_{\hbox{mre}})=R(\mu, {\hat q}_{\hbox{mre}})\equiv R_0$,
it is sufficient to show that
$\liminf_{k\to\infty}r_k(\pi_k,{\hat q}_k^\pi)\geq R_0$.
Making the transformations $z=x-\mu$ ($dz=dx$), $t=a-\mu$ ($d t=d a$) and $\xi_i=\mu_i/k$ ($d\xi_i=d\mu_i/k$) gives that
\begin{align*}
r_k(\pi_k, {\hat q}_k^\pi) =&
\int_{|\xi_i|<1/2, i=1,\ldots,p} \int\int \{ q(y) - {\hat q}_k^{\pi*}(y|x)\}^2 \;dy \,p(x) dx \, d\mu
\\
\geq&
\int_{|\xi_i|<(1-\epsilon)/2, i=1,\ldots,p} \int\int \{ q(y) - {\hat q}_k^{\pi*}(y|x)\}^2 \;dy \, p(x) \, dx\, d\mu,
\end{align*}
for any $\epsilon>0$, where
$$
{\hat q}_k^{\pi*}(y|x) = \int_{t+k\xi\in A_k} q(y-t)\, p(x-t)\, dt/\int_{t+k\xi\in A_k} p(x-t)\,dt.
$$
For $|\xi_i|<(1-\epsilon)/2$, it is seen that $\{t+k\xi\in A_k\} \supset \{-k\epsilon/2 < t_i <k\epsilon/2, i=1, \ldots, p\}$, which implies that ${\hat q}_k^{\pi*}(y|x)\to {\hat q}_{mre}(y|x)$ as $k\to\infty$.
Using Fatou's lemma, one gets
\begin{align*}
\liminf_{k\to\infty} r_k(\pi_k, {\hat q}_k^\pi) 
\geq&
\, \liminf_{k\to\infty} \int_{|\xi_i|<(1-\epsilon)/2, i=1,\ldots,p} \int\int \{ q(y) - {\hat q}_k^{\pi*}(y|x)\}^2 \;dy \, p(x)\, dx\,d\mu,
\\
\geq& 
\int_{|\xi_i|<(1-\epsilon)/2, i=1,\ldots,p} \int\int \liminf_{k\to\infty} \{ q(y) - {\hat q}_k^{\pi*}(y|x)\}^2 \;dy \, p(x)\, dx\,d\mu,
\\
=& \;
(1-\epsilon)^p R(\mu, {\hat q}_{mre}(y|x))= (1-\epsilon)^p R_0\,.
\end{align*}
From the arbitrariness of $\epsilon>0$, it follows that 
$\liminf_{k\to\infty} r_k(\pi_k, {\hat q}_k^\pi) \geq R_0$, which proves the minimaxity of the MRE predictor.  \qed

\section*{Acknowledgements}

Tatsuya Kubokawa's research is supported in part by Grant-in-Aid for Scientific Research Nos. 19200020 and 21540114 from the Japan Society for the Promotion of Science, 
Eric Marchand's research is supported in part by the Natural Sciences and Engineering Research Council of Canada, and William Strawderman's research is partially supported by a grant from the Simons Foundation (\#209035).  Finally, thanks to the Center for International Research on the Japanese Economy which provided financial support for a 2012 visit by Marchand to the University of Tokyo.

\renewcommand{\baselinestretch}{0.8}

\end{document}